\documentclass[11pt, reqno]{amsart}

\usepackage{hyperref,verbatim,xcolor,graphicx,float,lmodern,enumerate,cite,enumitem,pgfplots}
\usepackage{amsmath,amsthm,mathtools,amssymb,bbm,subfig,wrapfig,sepfootnotes}
\usepackage{graphicx}
\graphicspath{ {./} }

\numberwithin{equation}{section}

\theoremstyle{plain}

\numberwithin{equation}{section}

\newtheorem{theorem}{Theorem}[section]
\newtheorem{lemma}[theorem]{Lemma}

\newtheorem{corollary}[theorem]{Corollary}
\newtheorem{remark}[theorem]{Remark}

\newtheorem{assumption}[theorem]{Assumption}

\newcommand{\dmin}{m_0}
\newcommand{\dmax}{m_\infty}
\newcommand{\Exp}{\mathbb{E}}
\newcommand{\Pro}{\mathbb{P}}
\newcommand{\lep}{\left(}
\newcommand{\rip}{\right)}
\newcommand{\eee}{\mathrm{e}}

\def\CM{\mathrm{CM}}

\newcommand{\R}{\mathbb{R}}
\newcommand{\Z}{\mathbb{Z}}
\newcommand{\N}{\mathbb{N}}

\makeatletter
\let\oldabs\abs
\def\abs{\@ifstar{\oldabs}{\oldabs*}}

\begin{document}
\title[Largest eigenvalue of configuration model]{Largest eigenvalue of the configuration model\\ 
and breaking of ensemble equivalence}

\author[P. Dionigi]{Pierfrancesco Dionigi}
\address{Mathematical Institute, Leiden University, P.O.\ Box 9512, 2300 RA Leiden, The Netherlands.}
\email{p.dionigi@math.leidenuniv.nl}
\author[D. Garlaschelli]{Diego Garlaschelli}
\address{Lorentz Institute for Theoretical Physics, Leiden University, P.O.\  Box 9504, 2300 RA Leiden, The Netherlands
\& IMT School for Advanced Studies, Piazza S.\ Francesco 19, 55100 Lucca, Italy
\& INdAM-GNAMPA Istituto Nazionale di Alta Matematica, Italy}
\email{garlaschelli@lorentz.leidenuniv.nl}
\author[R.~S.~Hazra]{Rajat Subhra Hazra}
\author[F.\ den Hollander]{Frank den Hollander}
\address{Mathematical Institute, Leiden University, P.O.\ Box 9512, 2300 RA Leiden, The Netherlands}
\email{r.s.hazra@math.leidenuniv.nl}
\email{denholla@math.leidenuniv.nl}

\begin{abstract}
We analyse the largest eigenvalue of the adjacency matrix of the configuration model with large degrees, where the latter are treated as hard constraints. In particular, we compute the expectation of the largest eigenvalue for degrees that diverge as the number of vertices $n$ tends to infinity, uniformly on a scale between $1$ and $\sqrt{n}$, and show that a weak law of large numbers holds. We compare with what was derived in our earlier paper~\cite{dionigiCentralLimitTheorem2023} for the Chung-Lu model, which in the regime considered represents the corresponding configuration model with soft constraints, and show that the expectation is shifted down by $1$ asymptotically. This shift is a signature of breaking of ensemble equivalence between the hard and soft (also known as micro-canonical and canonical) versions of the configuration model. The latter result generalizes a previous finding~\cite{dionigiSpectralSignatureBreaking2021} obtained in the case when all degrees are equal.
\end{abstract}
\thanks{PD, RSH, FdH and MM are supported by the Netherlands Organisation for Scientific Research (NWO) through the Gravitation-grant NETWORKS-024.002.003. DG is supported by the Dutch Econophysics Foundation (Stichting Econophysics, Leiden, The Netherlands) and by the European Union - NextGenerationEU - National Recovery and Resilience Plan (Piano Nazionale di Ripresa e Resilienza, PNRR), project `SoBigData.it - Strengthening the Italian RI for Social Mining and Big Data Analytics' - Grant IR0000013 (n. 3264, 28/12/2021) (\url{https://pnrr.sobigdata.it/}).}

\date{\today}

\maketitle

\section{Introduction and main results}


\paragraph{\bf Motivation.}

Spectral properties of adjacency matrices in random graphs play a crucial role in various areas of network science. The largest eigenvalue is particularly sensitive to the graph architecture, making it a key focus. In this paper we focus on a random graph with a hard constraint on the degrees of nodes. In the homogeneous case (all degrees equal to $d$), it reduces to a random $d$-regular graph. In the heterogeneous case (different degrees), it is known as the configuration model. Our interest is characterizing the expected largest eigenvalue of the configuration model and comparing it with the same quantity for a corresponding random graph model where the degrees are treated as soft constraints.

The set of $d$-regular graphs on $n$ vertices, with $d=d(n)$, is non-empty when $1 \leq d \leq n-1$ and $dn$ is even. Selecting a graph uniformly at random from this set results in what is known as the \emph{random $d$-regular graph}, denoted by $G_{n,d}$. The spectral properties of $G_{n,d}$ are well-studied for $d\ge 2$ (for $d=1$, the graph is trivially a set of disconnected edges). For instance, all eigenvalues of its adjacency matrix fall in the interval $[-d,d]$, with the largest eigenvalue being $d$. The computation of $\lambda = \max\{|\lambda_2|, |\lambda_n|\}$, where $\lambda_2$ and $\lambda_n$ are the second-largest and the smallest eigenvalue, respectively, has been challenging. It is well known that for fixed $d$ the empirical distribution function of the eigenvalues of $G_{n,d}$ converges to the so-called Kesten-McKay law~\cite{mckay1981expected}, the density of which is given by
\begin{equation*}
f_{\rm KM}(x) = 
\begin{cases} 
\frac{d\sqrt{4(d-1)-x^2}}{2\pi (d^2-x^2)} &\quad \text{if } |x| \le 2\sqrt{d-1},\\
0 &\quad \text{otherwise}.
\end{cases}
\end{equation*}
A lower bound of the type $\lambda \ge 2\sqrt{d-1}$ was expected from this, and was indeed proven in the celebrated Alon-Boppana theorem in \cite{alonEigenvaluesExpanders1986}, stating that $\lambda \ge 2\sqrt{d-1}+O_d(\log^{-2}(n))$, where the constant in the error term only depends on $d$. In the same paper an upper bound of the type $\lambda \leq 2\sqrt{d-1}+o(1)$ was conjectured . This conjecture was later proven in the pioneering work \cite{Friedman} (for a shorter proof, see \cite{Bordenave:newproof}). When $d=d(n)\to\infty$, the spectral distribution converges to the semi-circular law \cite{tran:vu:wang}. It was proven in \cite{broderOptimalConstructionEdgeDisjoint1998} that $\lambda = O(\sqrt{d}\,)$ with high probability when $d=o(\sqrt{n}\,)$. This result was extended in \cite{BHKY} to $d=o(n^{2/3})$. For recent results and open problems, we refer the reader to \cite{Vu:survey}.

There has not been much work on the inhomogeneous setting where not all the degrees of the graph are the same. The natural extension of the regular random graph is the configuration model, where a degree sequence $\vec{d}_n=(d_1,\ldots,d_n)$ is prescribed. Unless otherwise specified, we assume the degrees to be hard constraints, i.e. realized exactly on each configuration of the model (this is sometimes called the `microcanonical' configuration model, as opposed to the `canonical' one where the degrees are soft constraints realized only as ensemble averages~\cite{squartini2015breaking,garlaschelliCovarianceStructureBreaking2018}). The empirical spectral distribution of the configuration model is known under some assumptions on the growth of the sum of the degrees. When $\sum_{i=1}^n d_i = O(n)$, the graph locally looks like a tree. It has been shown that the empirical spectral distribution exists and that the limiting measure is a functional of a size-biased Galton-Watson tree \cite{Bordenave:LeLarge}. When $\sum_{i=1}^n d_i$ grows polynomially with $n$, the local geometry is no longer a tree. This case was recently studied in \cite{demboEmpiricalSpectralDistributions2021}, where it is shown that the appropriately scaled empirical spectral distribution of the configuration model converges to the free multiplicative convolution of the semi-circle law with a measure that depends on the empirical degree sequence. The scaling of the largest and the second-largest eigenvalues has not yet been studied in full generality.

For the configuration model the behaviour of the largest eigenvalue is non-trivial. In the present paper, we consider the configuration model with large degrees, compute the expectation of the largest eigenvalue of its adjacency matrix, and prove a weak law of large numbers. When $d\to\infty$, the empirical distribution of $G_{n,d}$ (appropriately scaled) converges to the semi-circle law. Also, if we look at Erd\H{o}s-R\'enyi random graphs on $n$ vertices with connection probability $d/n$, then the appropriately scaled empirical distribution converges to the semi-circle law. It is known that the latter two random graphs exhibit different behaviour for the expected largest eigenvalue \cite{dionigiSpectralSignatureBreaking2021}. This can be understood in the broader perspective of \emph{breaking ensemble equivalence}~\cite{squartini2015breaking,garlaschelliCovarianceStructureBreaking2018}, which we discuss later. In the inhomogeneous setting, the natural graph to compare the configuration model with is the Chung-Lu model. For the case where the sum of the degrees grows with $n$, it was shown in \cite{Chakrabarty:Hazra:denHollander:Sfragara} that the empirical spectral distribution converges to the free multiplicative product of the semi-circle law with a measure that depends on the degree sequence, similar to \cite{demboEmpiricalSpectralDistributions2021}. In \cite{dionigiCentralLimitTheorem2023}, we investigated the largest eigenvalue and its expectation for the Chung-Lu model (and even derived a central limit theorem). In the present paper we complete the picture by showing that, when the degrees are large, there is a \emph{gap} between the expected largest eigenvalue of the Chung-Lu model and the configuration model. We refer the reader to the references listed in \cite{dionigiSpectralSignatureBreaking2021,dionigiCentralLimitTheorem2023} for more background.


\paragraph{\bf Main theorem.}

For $n\in\N$, let $\CM(\vec{d}_n)$ be the random graph on $n$ vertices generated according to the \emph{configuration model} with degree sequence $\vec{d}_n= (d_i)_{i \in [n]} \in \mathbb{N}^n$ \cite[Chapter 7.2]{vanderhofstadRandomGraphsComplex2016}. Define
\begin{equation*}
m_0= \min_{i\in [n]} d_i, \qquad m_\infty = \max_{i \in [n]} d_i, \qquad m_k = \sum_{i \in [n]} (d_i)^k, \quad k \in \N.
\end{equation*}
Throughout the paper we need the following assumptions on $\vec{d}_n$ as $n\to\infty$.

\begin{assumption}
\label{ass:deg}
$\mbox{}$\\
{\rm 
(D1) {\bf Bounded inhomogeneity:} $m_0 \asymp m_\infty$.\\
(D2) {\bf Connectivity and sparsity:} $1 \ll m_\infty \ll \sqrt{n}$.
}
\hfill$\spadesuit$
\end{assumption}

\noindent
Under these assumptions, $\CM(\vec{d}_n)$ is with high probability \emph{non-simple} \cite[Chapter 7.4]{vanderhofstadRandomGraphsComplex2016}. We write $\Pro$ and $\Exp$ to denote probability and expectation with respect to the law of $\CM(\vec{d}_n)$ \emph{conditional on being simple}, suppressing the dependence on the underlying parameters. 

Let $A_{\CM(\vec{d}_n)}$ be the adjacency matrix of $\CM(\vec{d}_n)$. Let $(\lambda_i)_{i\in[n]}$ be the eigenvalues of $A_{\CM(\vec{d}_n)}$, ordered such that $\lambda_1 \geq\dots\geq \lambda_n$. We are interested in the behaviour of $\lambda_1$ as $n\to\infty$. Our main theorem reads as follows.

\begin{theorem}
\label{thm:main}
Subject to Assumption~\ref{ass:deg},
\begin{equation}
\label{eq:mic}
\Exp\left[\lambda_1\right] 
= \frac{m_2}{m_1} + \frac{m_1m_3}{m_2^2} -1 + o(1), \qquad n \to \infty,
\end{equation}
and $$\frac{\lambda_1}{\Exp\left[\lambda_1\right]} \to 1 \text{ in $\Pro$-probability. }$$
\end{theorem}

In \cite{dionigiCentralLimitTheorem2023} we looked at an alternative version of the configuration model, called the \emph{Chung-Lu model} $\CM^*_n(\vec{d}_n)$, where the \emph{average degrees}, rather than the degrees themselves, are fixed at $\vec{d}_n$. This is an ensemble with soft constraints; in the considered regime for the degrees, it coincides with a maximum-entropy ensemble, also called `canonical' configuration model~\cite{squartini2015breaking}. For this model we showed that, subject to Assumption~\ref{ass:deg},
\begin{equation}
\label{eq:can}
\Exp^*\left[\lambda_1\right] 
= \frac{m_2}{m_1} + \frac{m_1m_3}{m_2^2} + o(1), \qquad n \to \infty,
\end{equation}
and $\lambda_1/\Exp^*\left[\lambda_1\right] \to 1$ in $\Pro^*$-probability, where $\Pro^*$ and $\Exp^*$ denote expectation with respect to the law of $\CM^*_n(\vec{d}_n)$ and $\lambda_1$ is the largest eigenvalue of the $A_{\CM^*_n(\vec{d}_n)}$. The notable difference between \eqref{eq:mic} and \eqref{eq:can} is the shift by $-1$. 

For the special case where all the degrees are equal to $d$, we have $m_0=m_\infty=d$ and $m_k = nd^k$, and so $\Exp[\lambda_1] = d+o(1)$ and $\Exp^*[\lambda_1] = d+1+o(1)$. In fact, $\Pro(\lambda_1=d)=1$. Since in this model the degrees can fluctuate with the same law (\emph{soft} constraint in the physics literature), this case reduces to the Erd\H{o}s-R\'enyi random graph for which results on $\Exp^*[\lambda_1] $ were already well known in \cite{juhasz1981spectrum,furedi1981eigenvalues} and further analyzed in \cite{erdosSpectralStatisticsErdos2013b}.


\paragraph{\bf Breaking of ensemble equivalence.}
 
The shift by $-1$ was proven in \cite{dionigiSpectralSignatureBreaking2021} for the homogeneous case with equal degrees and is a \emph{spectral signature} of \emph{breaking of ensemble equivalence}~\cite{squartini2015breaking}. Indeed, a $d$-regular graph is the `micro-canonical' version of a random graph where all degrees are equal and `hard', and the Erd\H{o}s-R\'enyi random graph is the corresponding `canonical' version where all degrees are equal and `soft'. 
More in general, $\CM(\vec{d}_n)$ is the micro-canonical configuration model where the constraint on the degrees is `hard', while $\CM^*_n(\vec{d}_n)$ is the canonical version where the constraint is `soft'. 
We refer the reader to \cite{garlaschelliCovarianceStructureBreaking2018} for the precise definition of these two configuration model \emph{ensembles} and for the proof that they are not asymptotically equivalent \emph{in the measure sense}~\cite{touchette2015equivalence}. This means that the \emph{relative entropy per node} of $\Pro$ with respect to $\Pro^*$ has a strictly positive limit as $n \to \infty$. This shows that the choice of constraint matters, not only on a microscopic scale but also on a macroscopic scale. Indeed, for non-equivalent ensembles one expects the existence of certain macroscopic properties that have different expectation values in the two ensembles (\emph{macrostate (in)equivalence}~\cite{touchette2015equivalence}). 
The fact that the largest eigenvalue picks up this discrepancy is interesting. What is remarkable is that the shift by $-1$, under the hypotheses considered, holds true also in the case of heterogeneous degrees and remains the same \emph{irrespective} of the scale of the degrees and of the distribution of the degrees on this scale.


\paragraph{\bf Outline.}

The remainder of this paper is organised as follows. In Section \ref{sec:simple} we look at the issue of simplicity of the graph. In Section \ref{sec:bounds} we bound the spectral norm of the matrix
\begin{equation}
\label{Hdef}
H = A_{\CM(\vec{d}_n)} - \Exp[A_{\CM(\vec{d}_n)}].
\end{equation}
We use the proof of \cite{broderOptimalConstructionEdgeDisjoint1998} to show that $\|H\|= \mathrm{o}(m_{\infty})$ with high probability. Using the latter we show that
\begin{equation*}
\lambda_1 \sim \lambda_1\left(\Exp[A_{\CM(\vec{d}_n)}]\right) \sim \frac{m_2}{m_1}
\end{equation*} 
with high probability. In Section \ref{sec:expansion} we use the estimates in Section \ref{sec:bounds} to prove Theorem \ref{thm:main}.


\section{Configuration model and simple graphs with a given degree sequence}
\label{sec:simple}

In the context of random graphs with hard constraints on the degrees, we work with both the \emph{configuration model} (a random \emph{multi-graph} with a prescribed degree sequence) and with the \emph{conditional configuration model} (a random \emph{simple} graph with a prescribed degree sequence). We view the second as a special case of the first.


\paragraph{\bf Configuration model.}
 
The configuration model with degree sequence $\vec{d}=(d_1,\dots,d_n)$, $\CM(\vec{d}_n)$, generates a graph through a \emph{perfect matching} $P$ of the set of half-edges $\mathcal{E} = \cup_{i \in [n]} \{i\} \times [d_i]$, i.e., $P\colon\, \mathcal{E} \to \mathcal{E}$ such that $P$ is an isomorphism and an involution, and $P(\alpha) \neq \alpha$ for all $\alpha \in \mathcal{E}$. 

\begin{figure}[htbp]
\centering
\subfloat[][\emph{Pairing procedure}.]
{\includegraphics[width=0.4\textwidth]{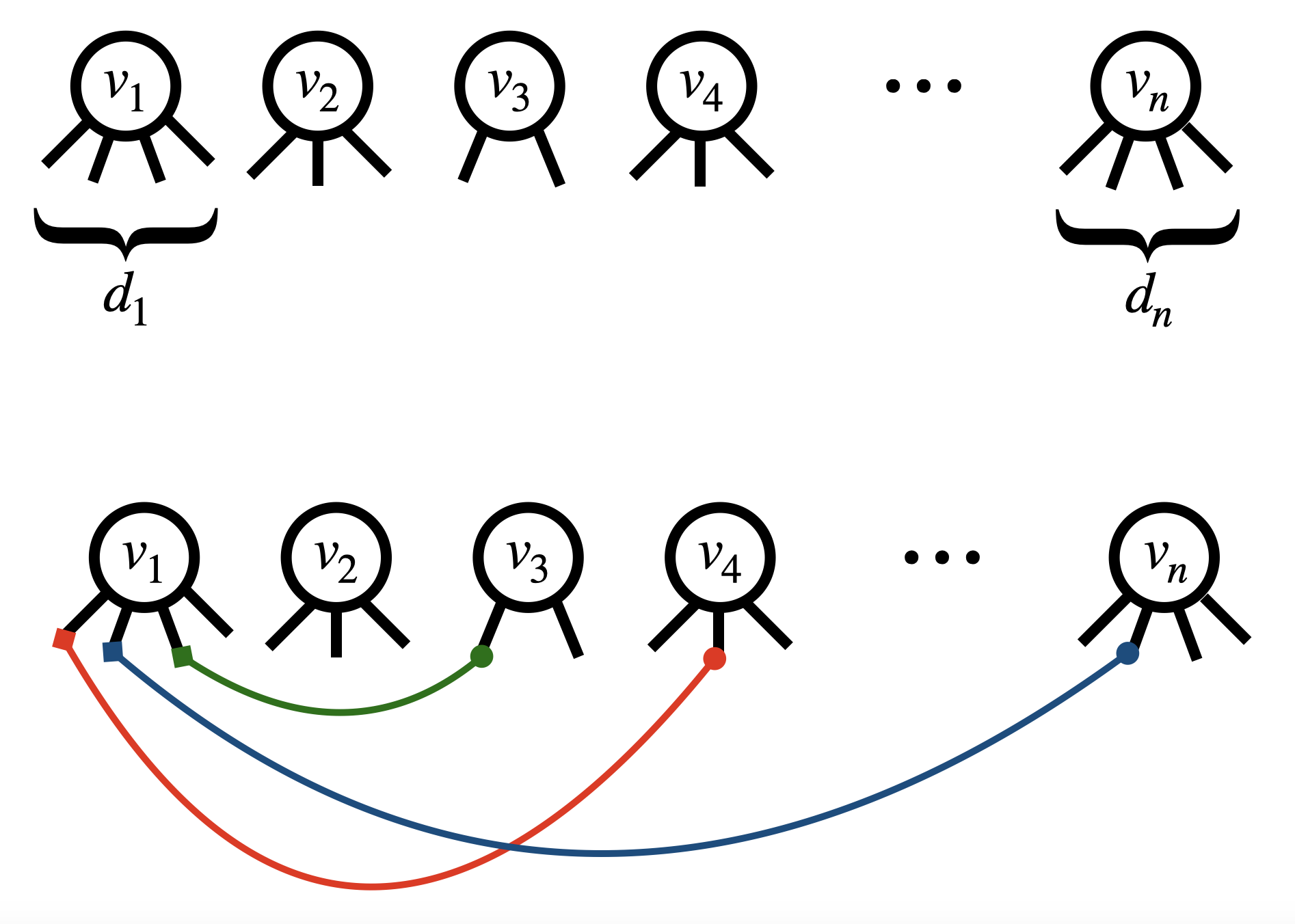}} \quad
\subfloat[][\emph{The configuration model generates multi-graphs with self-loops (left) and multiple edges (right)}.]
{\includegraphics[width=0.4\textwidth]{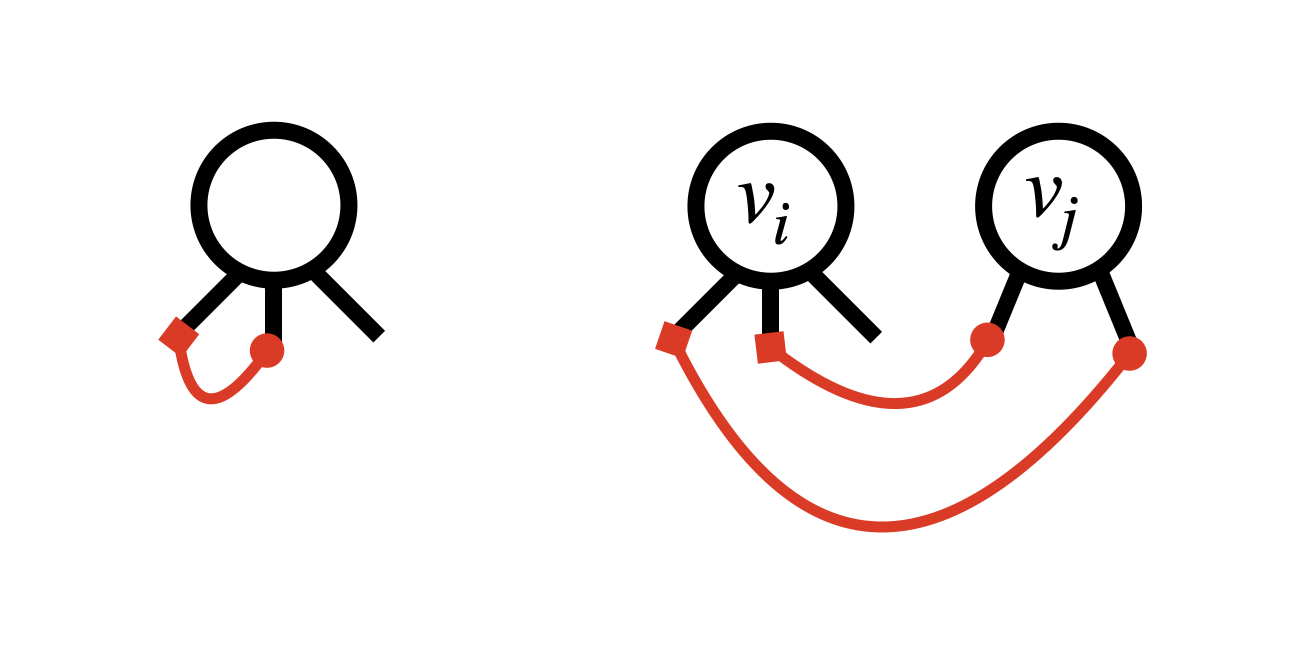}}
\label{fig:confmodel}
\end{figure}

It turns out that a pairing scheme, such as matching half-edges from left to right or selecting pairs uniformly at random, is considered admissible as long as it yields the correct probability for a perfect matching (see, for example, \cite[Lemma 7.6]{vanderhofstadRandomGraphsComplex2016}). An important property of each admissible pairing scheme is that
\begin{equation}
\label{eq:pairscheme}
\Pro_{\CM(\vec{d}_n)}\lep\alpha\sim\beta\rip=\frac{1}{m_1-1}\qquad \forall \alpha,\beta\in\mathcal{E}.
\end{equation}
We can therefore endow $\mathcal{E}$ with the lexicographic order. An element $\alpha \in \mathcal{E}$, $\alpha=(i,\ell)$, is associated with a vertex $v(\alpha)=i$ through the first component of $\alpha$. An edge $e$ is an element of $\mathcal{E}\times \mathcal{E}$ given by the pairing $e = \{\alpha,\beta\} = \{\alpha,P(\alpha)\}$. We define the configuration $\mathcal{C}$ as the set of all such $e$. Given the lexicographic order, we may assume that $\alpha < \beta$ and identify arbitrarily the head and the tail of the edge: $t(e)=v(\alpha)$, $h(e)=v(\beta)$. We can then order the edges in their order of appearance via $t(e)$, forming a list $\{e_i, \, i=1, \ldots, m_1/2\}$. These properties of the configuration model will be used in Subsection \ref{sec:spec}. For the configuration model it is easy to check that\footnote{We adopt the convention that a self-loop contributes $2$ to the degree, i.e., $a_{ii}$ is twice the number of self-loops attached to vertex $i$. This convention is useful because it yields $\sum_j a_{ij} = d_i$.} 
\begin{equation}
\label{eq:expectationCM}
\Exp[a_{ii}] = \frac{d_i(d_i-1)}{m_1-1}, \qquad \Exp[a_{ij}] = \frac{d_id_j}{m_1-1}, \quad i \neq j. 
\end{equation}
Indeed every matrix element $a_{ij}$, can be expressed as 
\begin{equation*}
\sum_{k=1}^{d_i}\sum^{d_j}_{h=1} \mathbbm{1}\big(\alpha\sim\beta,\alpha=(i,k),\beta=(j,h)\big),
\end{equation*}
where $\mathbbm{1}(E)$ is the indicator function of the event $ E $. Taking expectations and using \eqref{eq:pairscheme}, we get \ref{eq:expectationCM}.


\paragraph{\bf Simple graphs with a given degree sequence.}

A linked but different model is the one that samples uniformly at random a simple graph with a given degree sequence $\vec{d}$, $\mathcal{G}(\vec{d})$\footnote{With a little abuse of notation we permit our simple graphs to have self-loops. We maintain the convention explained in the previous footnote.}. An immediate question is whether there is any relation between $\CM(\vec{d})$ and $\mathcal{G}(\vec{d})$. It turns out that the following is true (see \cite{vanderhofstadRandomGraphsComplex2016} for reference):
\begin{lemma}
Let be $G_n$ a graph generated via the configuration model $\CM(\vec{d}_n)$ with degree sequence $\vec{d}$. Then
\begin{equation}
\label{lem:simpleCM_1}
\Pro_{\CM(\vec{d}_n)}(G_n \textnormal{ is simple}) = \exp\left[-O\left(\frac{m^2_1}{n^2}\right)\right]
\end{equation}
and
\begin{equation}
\label{lem:simpleCM_2}
\Pro_{\CM(\vec{d}_n)}\lep G_n = \cdot \mid G_n\textnormal{ is simple}\rip \textnormal{ is uniform}.
\end{equation} 
\end{lemma}

\noindent
Thus, we can identify the law of the \emph{uniform random graph with a given degree sequence} with the law of the \emph{configuration model conditional on simplicity}, i.e., 
\begin{equation*}
\Pro_{\mathcal{G}(\vec{d})}(\cdot) = \Pro_{\CM(\vec{d})}(\cdot \mid \text{G is simple}).
\end{equation*} 
It follows that the expectation under the uniform random graph with a given degree sequence can be expressed as a the expectation of the configuration model conditional on simplicity. As explained in \cite[Remark 1.14]{vanderhofstadRandomGraphsComplex2023}, we have
\begin{equation}
\begin{split}
\label{eq:expectationUG}
\Exp_{\mathcal{G}(\vec{d})}[a_{ij}] 
&= \Pro_{\mathcal{G}(\vec{d})}(i\sim j) = (1+o(1)) \frac{d_id_j}{m_1 +d_i d_j}\\
&\leq \frac{d_id_j}{m_1(1+O(1/m_\infty))}=\frac{d_id_j}{m_1}+O(1/m_\infty). 
\end{split}
\end{equation}


\section{Bound on $\|H\|$}
\label{sec:bounds}

In this section we derive the following bound on the spectral norm of the matrix $H$ defined in \eqref{Hdef}. This bound will play a crucial role in the proof of Theorem~\ref{thm:main} in Section~\ref{sec:expansion}. 

\begin{theorem}
\label{thm:edge}
For every $K>0$, there exists a constant $C>0$ such that
\begin{equation*}
\|H\| \le C \sqrt{m_\infty}\, \text{ with probability } 1-o(n^{-K}),
\end{equation*} 
where the constant $C$ depends  on $\frac{m_\infty}{m_0}$ and $K$ only.
\end{theorem}

\noindent

The proof is given in Section~\ref{sec:spec} and we follow the proof of Lemma 16 of \cite{broderOptimalConstructionEdgeDisjoint1998}. In Sections~\ref{sec:first}--\ref{sec:second} we derive some estimates that are needed in Section~\ref{sec:spec}.

We have
\begin{equation*}
\|H\| = \max\{\lambda_1(H),\left|\lambda_n(H)\right|\}.
\end{equation*}
Recall \eqref{Hdef}. In Lemma \ref{lem:rank1} below we will see that $ \Exp[A_{\CM(\vec{d}_n)}]$ is asymptotically rank $1$ (more precisely, the other eigenvalues are of order $O(n^{-1})$. This means that $\lambda_i\lep A_{\CM(\vec{d}_n)}\rip = O(n^{-1})$, $i \neq 1$. Therefore, by Weyl's inequality, it easy to see that $\lambda_1(H)$ and $\lambda_n(H)$ are bounded above and below by $\lambda_2(A_{G_n})$ and $\lambda_n(A_{G_n})$ with error $O(n^{-1})$, while $\lambda_1(A_{G_n})$ is in a window of size $\sqrt{m_\infty}$ around $\frac{m_2}{m_1-1}$. We therefore have the following Corollary, which will be needed in Section \ref{sec:expansion}:

\begin{corollary}\label{cor:lambdaconc}
\label{cor:concentration}
For every $K>0$,
\begin{equation*}
\lambda_1(A_{G_n}) - \lambda_1(\Exp [A_{G_n}]) = O(\sqrt{m_\infty}\,) \text{ with probability } 1-o(n^{-K}). 
\end{equation*}
\end{corollary}

\subsection{Spectral estimates}
\label{sec:spec}


\paragraph{\bf Notation.}

Abbreviate $G_n = \CM(\vec{d}_n)$. Let $U$ be uniformly distributed on $[n]$, let $d_U$ be the degree of a vertex that is picked uniformly at random, and let
\begin{equation*}
\omega_n=\Exp[d_U]
\end{equation*}
be the average of the empirical degree distribution. Under Assumption \ref{ass:deg}, $\omega_n\to \infty$ and $\omega_n=o(n)$. We define the \emph{normalised degree sequence} $(\hat{d}_i)_{i\in[n]}$ as 
\begin{equation}
\label{eq:normalized}
\hat{d}_i=\frac{d_i}{\omega_n}
\end{equation}
and the \emph{normalised adjacency matrix} $\hat{A}_{G_n}$ as 
\begin{equation*}
\hat{A}_{G_n} = \frac{A_{G_n}}{\sqrt{\omega_n}}.
\end{equation*}
In the following we will need multiples of the vector
\begin{equation}
\label{eq:deftildee}
\tilde{e}_i=\frac{d_i}{\sqrt{m_1-1}}, \,\, 1\le i \le n,
\end{equation}
which is the eigenvector corresponding to the rank-$1$ approximation of $\Exp[A_{G_n}]$, as it easy to check using \eqref{eq:expectationCM}.


\paragraph{\bf Two lemmas.}

The matrix $\Exp[A_{G_n}]$ is asymptotically rank 1:

\begin{lemma}\label{lem:rank1}
\label{lem:shift}
Define 
\begin{equation*}
\tilde{e} = \lep\tilde{e}_1,\dots,\tilde{e}_n\rip^t
\end{equation*} 
and 
\begin{equation*}
A_\textnormal{shift} = \Exp[A_{G_n}]-|\tilde{e}\rangle\langle\tilde{e}|
= - \frac{1}{m_1-1}\,\textnormal{diag}\lep d_1,\dots,d_n\rip.
\end{equation*}
Then
\begin{equation*}
\|A_\textnormal{shift}\| = O(n^{-1}).
\end{equation*}
\end{lemma}

\begin{proof}
Using \eqref{eq:expectationCM}, we have $A_\textnormal{shift}$ is a diagonal matrix and hence $\|A_\textnormal{shift}\|\leq \frac{m_\infty}{m_1-1} = O(n^{-1})$ under Assumption~\ref{ass:deg}. Note that when all the degrees are equal to $d$, then this bound is sharp, i.e., $\|A_\textnormal{shift}\|=\frac{d}{dn-1}$, and $\Exp [A_{G_n}]$ is exactly rank 1.
\end{proof}

It is shown in \cite{demboEmpiricalSpectralDistributions2021} that the empirical spectral distribution of $A_{G_n}$ has a deterministic limit given by $\mu=\mu_{sc} \boxtimes \mu_{\hat{D}}$, where $\mu_{sc}$ is the standard Wigner semicircle law, $\mu_{\hat{D}}$ is the distribution of $\hat{d}_U$, and $\boxtimes$ is the free product defined in \cite{voiculescuMultiplicationCertainNonCommuting1987,bercoviciFreeConvolutionMeasures1993,arizmendie.TransformSymmetricProbability2009}. 

\begin{lemma}
Let $(\hat{d}_i)_{i \in [n]}$ be the normalised degree sequence, and suppose that Assumption~\ref{ass:deg} holds and
\begin{equation}
\label{eq:deg_conv}
\frac{1}{n}\sum_{i=1}^n \delta_{\hat{d}_i}\Rightarrow \mu_{\hat D}.
\end{equation}
Then $\mu_{\hat D}$ is compactly supported.
\end{lemma}

\begin{proof}
By Assumption \ref{ass:deg} D(2),
\begin{equation*}
0 < c \leq \liminf_{n\to\infty} \frac{\min_i\hat{d}_i}{\max_i\hat{d}_i} \leq 
\limsup_{n\to\infty} \frac{\max_i\hat{d}_i}{\min_i\hat{d}_i} \leq C < \infty.
\end{equation*}
Hence the support of $\hat d_U$ is contained in a multiple of $[c,C]$, and \eqref{eq:deg_conv} implies that $\mu_{\hat D}$ is compactly supported.
\end{proof}


\paragraph{\bf Key estimates.}

From \cite[Theorem 1]{demboEmpiricalSpectralDistributions2021} we know that, under Assumption~\ref{ass:deg}, the empirical spectral distribution of $\hat{A}_{G_n}$, written
\begin{equation*}
\mu_{\hat{A}_{G_n}},
\end{equation*} 
converges weakly to $\mu_{\text{sc}} \boxtimes \mu_{\hat{D}}$. The fact that $\mu_{sc}$ and $\mu_{\hat{D}}$ are compactly supported implies that $\mu_{\text{sc}} \boxtimes \mu_{\hat{D}}$ is compactly supported too. However, this still allows that $H$ has $o(n)$ outliers, which possibly control $\|H\|$. It is hard to analyse $\|H\|$ for sparse graphs, because it is related to expansion properties of the graph, mixing times of random walks on the graph, and more. To prove that $\|H\|= O(\sqrt{m_\infty}\,)$, we will use the argument described in \cite{broderOptimalConstructionEdgeDisjoint1998}, which is an adaptation of the argument in \cite{10.1145/73007.73063}. We are only after the order of $\|H\|$, not sharp estimates. 

For random regular graphs with fixed degree, the problem of computing $\|H\|$ was solved in \cite{friedmanProofAlonSecond2004}. We are not aware of any proof concerning the second largest eigenvalue of configuration models with bounded degrees, but any sharp bound must come from techniques of the type employed in \cite{friedmanProofAlonSecond2004,bordenaveNewProofFriedman2019}. In what follows we prove Theorem~\ref{thm:edge} based on \cite{10.1145/73007.73063} and its adaptation to the inhomogeneous setting in \cite{broderOptimalConstructionEdgeDisjoint1998}. In the following, we present the proof as outlined in the latter paper, adapting it to our case to make our paper self-contained. This adjustment is necessary because the result we need is somewhat obscured in the context in which it appears in \cite{broderOptimalConstructionEdgeDisjoint1998}.

Let $D_n=\text{diag}(d_1,\dots, d_n)$. The transition kernel of a random walk on the graph $G_n$ is given by $P_n=D_n^{-1}A_{G_n}$. The \emph{random} matrix $P_n$ has as principal normalised eigenvector $\vec{1}=(1,\dots,1)/\sqrt{n}$ with eigenvalue 1. Define 
\begin{equation*}
\lambda^* = \max\{\lambda_2(P_n), |\lambda_n(P_n)|\}.
\end{equation*}
Note that the matrices $P_n$ and $S_n=D_n^{1/2}P_nD_n^{-1/2}$ have the same spectrum by a similarity transformation. Hence we can write the Rayleigh formula
\begin{equation*}
\lambda^* = \max_{z\colon\,\langle z,\vec{1} \rangle = 0} \frac{|\langle P_nz,z \rangle|}{\|z\|^2}
= \max_{x\colon\,\langle x,\sqrt{\vec{d}} \rangle=0} \frac{|\langle S_nx,x \rangle|}{\|x\|^2},
\end{equation*}
where $\sqrt{\vec{d}}=( \sqrt{d_1}, \ldots, \sqrt{d_n})$. 
Define $M_n=D_n^{-1}A_nD_n^{-1}$. Let $\tilde \lambda$ be the second largest eigenvalue of $M_n$. Then
\begin{equation*}
\langle S_nx,x \rangle = \langle D_n^{1/2} M_n D_n^{1/2}x,x \rangle = \langle M_nD_n^{1/2}x, D_n^{1/2} x \rangle.
\end{equation*}
Since 
\begin{equation*}
\frac{\langle x,x \rangle}{\|x\|^2} \geq \frac{1}{m_\infty} \frac{\langle D_n^{1/2}x, D_n^{1/2}x \rangle}{\|x\|^2},
\end{equation*}
putting $y=D_n^{1/2}x$ we see that 
\begin{equation}
\label{eq:max_1}
\lambda^*\leq m_\infty \max_{\langle y, \vec{1} \rangle=0} \frac{|\langle M_ny,y \rangle |}{\|y\|^2}=m_\infty\tilde{\lambda},
\end{equation}
which gives a bound of the type $\|H\| = O(m_\infty^2\tilde{\lambda})$. Note that the matrix elements of $M_n$ can be expressed as
\begin{equation}
\label{eq:defM_n}
(M_n)_{ij}= \frac{a_{ij}}{d_i d_j},
\end{equation}
where $a_{ij}$ counts the number of edges between vertices $i$ and $j$, with the convention for the diagonal elements stated earlier. In view of \eqref{eq:max_1}, in order to obtain a bound on $\lambda^*$ we must focus on $\tilde{\lambda}$. In fact, we must show that 
\begin{equation}
\label{claim:main}
\tilde{\lambda} = O(m_\infty^{-3/2}),
\end{equation}
in order to obtain the desired bound $\|H\| = O(\sqrt{m_\infty}\,)$. To achieve this we proceed in steps: 
\begin{enumerate}
\item[--] 
We reduce the computation of $\tilde{\lambda}$ to the analysis of two terms.
\item[--]
In \eqref{eq:expectation_X} below we identify the leading order of $\Exp[\tilde{\lambda}]$ from these two terms, which turns out to be $O(m_\infty^{-3/2})$.
\item[--] 
We show that the other terms are of higher order and therefore are negligible. 
\item[--] 
We prove concentration around the mean through concentration of the leading order term in Lemma \ref{lem:concentration} below.
\item[--] 
Recalling the denominator of \eqref{eq:defM_n}, we get a bound on $\|H\|$ after multiplying by $m_\infty^2$. 
\end{enumerate} 

To prove \eqref{claim:main},  we reduce the problem to a maximisation problem in a simpler space. Namely, let $\varepsilon \in (0,1)$, and 
\begin{equation*}
T = \left\{x\in \lep \frac{\varepsilon}{\sqrt{n}}\mathbb{Z}\rip^n\colon\,\sum_{i\in [n]} x_i = 0,\,\sum_{i \in [n]} x_i^2\leq 1\right\}.
\end{equation*}
Then, using the formula for the volume of the $n$-dimensional ball, we have
\begin{equation*}
|T| \leq \lep\frac{(2+\varepsilon)\sqrt{n}}{2\varepsilon}\rip^n \frac{\pi^{n/2}}
{\Gamma(\tfrac{n}{2}+1)}\leq \lep\frac{(2+\varepsilon)\sqrt{2\pi e}}{2\varepsilon}\rip^n.
\end{equation*}
The maximisation over $\mathbb{R}^n$ can be reduced to over $T$ and this leads to an error that depends on $\varepsilon$.

\begin{lemma}
Let 
\begin{equation*}
\lambda = \max \{|\langle x,M_ny \rangle|\colon\, x,y\in T\}.
\end{equation*}
Then
\begin{equation*}
\tilde{\lambda} \leq (1-\varepsilon)^{-2} \lambda.
\end{equation*}
\end{lemma}

\begin{proof}
Let be $\mathcal{S} = \{x\in\mathbb{R}^n\colon\,\sum_{i \in [n]} x_i=0, \|x\|\leq1 \}$. We want to show that for every $x\in \mathcal{S}$ there is a vector $y\in T$ such that $\|x-y\|\leq\varepsilon$ and $\sum_{i \in [n]} (x_i-y_i)=0$. Let us write the components of $x$ as 
\begin{equation*}
x_i = \varepsilon \frac{m_i}{\sqrt{n}} + f_i, \qquad i \in [n],
\end{equation*}
where $m_i \in \Z$ and $f_i\in [0,\varepsilon n^{-1/2} )$ is an error term. Because $\sum_{i \in [n]} x_i=0$, we choose $m_i$ such that $\sum_{i \in [n]} f_iv=v\varepsilon f n^{-1/2}$, where $f$ is a non-negative integer smaller than $n$. We relabel the indices in such a way that $m_i \leq m_j$ when $i \leq j$. Now consider the vector $y$ given by 
\begin{equation*}
y_i = \left\{\begin{array}{ll}
\varepsilon \frac{m_i +1}{\sqrt{n}} &\text{if } i \leq f,\\
\varepsilon \frac{m_i}{\sqrt{n}}      &\text{if } i > f. 
\end{array}
\right.
\end{equation*}
It follows that $\sum_{i \in [n]} y_i = 0$ and $\|y\| \leq 1$, and therefore $y\in T$. Furthermore, because $|x_i-y_i| \leq \varepsilon n^{-1/2}$ by construction, we have the required property $\|x-y\| \leq 1$. Iterating the previous argument, we can express every vector $x\in\mathcal{S}$ in terms of a sequence of vectors $(y^{(i)})_{i \in [n]}$ in $T$ such that
\begin{equation*}
x = \sum_{i \in [n]} \varepsilon^i y^{(i)}.
\end{equation*}
Therefore, because \eqref{eq:max_1} is maximized on $\mathcal{S}$, we have 
\begin{equation*}
\langle x,M_nx \rangle = \sum_{i,j \in [n] \times [n]} \varepsilon^{i+j} \langle x^{(i)} M_n x^{(j)} \rangle
\leq \frac{1}{(1-\varepsilon)^2} \max\{ |\langle y,M_n z \rangle| \colon\,y,z \in T\},
\end{equation*}
from which the claim follows.
\end{proof}

Our next goal is to show that $|\langle x, M_n y\rangle |=o(m_\infty^{-3/2})$ for all $x,y \in T$ with a suitably high probability. This can be done in the following way. Fix $x,y\in T$ and define the random variable 
\begin{equation*}
X = \sum_{i,j \in [n] \times [n]} x_i (M_n)_{ij} y_j.
\end{equation*} 
Define the set of indices
\begin{equation*}
\mathcal{B} = \left\{(i,j) \in [n] \times [n]\colon\, 0 < |x_iy_j| < \frac{\sqrt{m_\infty}}{n}\right\}.
\end{equation*}
Then we can rewrite $X=X'+X''$ with
\begin{equation*}
X' = \sum_{(i,j) \in \mathcal{B}} x_i (M_n)_{ij} y_j,  \qquad X''=\sum_{(i,j) \notin \mathcal{B}} x_i (M_n)_{ij} y_j.
\end{equation*}
In Section~\ref{sec:first} we show that $\Exp[X']$ is of the correct order and that $X'$ is well concentrated around its mean. In Section~\ref{sec:second} we analyse $X''$, which is of a different nature and requires that we exclude subgraphs in the configuration model that are too dense.


\subsection{Estimate of first contribution}
\label{sec:first}

{\bf 1.}
Use \eqref{eq:expectationCM} and \eqref{eq:defM_n} to write out $\Exp[m_{ij}] =\Exp[(M_n)_{ij}]$. This gives
\begin{equation*}
\Exp[X'] = \sum_{(i,j) \in\mathcal{B}} \frac{x_i y_j}{m_1-1} -\sum_{(i,i) \in \mathcal{B}} \frac{x_i y_i }{d_i^2 (m_1-1)}.
\end{equation*}
In view of the bound on $|x_iy_i|$, the last term gives a contribution 
\begin{equation*}
\left|\sum_{(i,i) \in\mathcal{B}} \frac{x_i y_i }{d_i^2 (m_1-1)}\right| = O\left(\frac{1}{m_1 m_\infty^{3/2}}\right).
\end{equation*}
Since $x,y\in T$, we have that $\sum_{i \in [n]} x_i=0$ and $\sum_{j \in [n]} y_j=0$, and therefore $\sum_{(i,j) \in [n] \times [n]} x_i y_j = 0$. Hence
\begin{equation*}
\left|\sum_{(i,j) \in \mathcal{B}} x_i y_j\right| = \left|\sum_{(i,j) \notin \mathcal{B}} x_i y_j\right|, 
\end{equation*}
where we can bound the right-hand side as
\begin {equation*}
\left|\sum_{(i,j) \notin \mathcal{B}} x_i y_j\right| \leq \sum_{(i,j)\colon |x_i y_j| \geq \frac{\sqrt{m_\infty}}{n}}\left|x_i y_j\right| 
\leq \sum_{(i,j)\colon |x_i y_j | \geq \frac{\sqrt{m_\infty}}{n}} \frac{x_i^2y_j^2}{| x_i y_j |} 
\leq \frac{n}{\sqrt{m_\infty}} \sum_{(i,j)} x^2_i y_j^2\leq \frac{n}{\sqrt{m_\infty}}.
\end{equation*}
We can therefore conclude that 
\begin{equation}
\label{eq:expectation_X}
\left|\Exp[X']\right| \leq \frac{n}{(m_1-1)\sqrt{m_\infty}} + O\left(\frac{1}{m_1 m_\infty^{3/2}}\right).
\end{equation} 

\medskip\noindent
{\bf 2.}
To prove that $X'$ is concentrated around its mean, we use an argument originally developed in \cite{broderSecondEigenvalueRandom1987,shamirSharpConcentrationChromatic1987} and used for the configuration model in \cite{10.1145/73007.73063,broderOptimalConstructionEdgeDisjoint1998}. This argument employs the \emph{martingale structure} of the configuration model \emph{conditional on partial pairings}. Define 
\begin{equation*}
\chi(x)=\begin{cases}
x, \qquad \text{if } |x|<\frac{\sqrt{m_\infty}}{n},\\
0, \qquad \text{otherwise.}
\end{cases}
\end{equation*}
We can then express $X'$ as 
\begin{equation*}
X' = \sum_{e\in\mathcal{C}} \frac{\chi(x_{t(e)}y_{h(e)})}{d_{t(e)}d_{h(e)}}
+ \sum_{e\in\mathcal{C}} \frac{\chi(x_{h(e)}y_{t(e)})}{d_{t(e)}d_{h(e)}} = X_a' +X_b'.
\end{equation*}
We divide the set of pairings $\mathcal{C}$ into three sets
\begin{equation*}
\begin{split}
&\mathcal{C}_1 
= \left\{e\in \mathcal{C}\colon\,|x_{t(e)}|>\frac{1}{\varepsilon \sqrt{n}}\right\},\\
&\mathcal{C}_2 
= \left\{e\in \mathcal{C}\colon\,|y_{t(e)}|>\frac{1}{\varepsilon \sqrt{n}},\,|x_{t(e)}|\leq\frac{1}{\varepsilon \sqrt{n}}\right\},\\
&\mathcal{C}_3 
= \left\{e\in \mathcal{C}\colon\,|y_{t(e)}|\leq\frac{1}{\varepsilon \sqrt{n}},\,|x_{t(e)}|\leq\frac{1}{\varepsilon \sqrt{n}}\right\},
\end{split}
\end{equation*}
and write
\begin{equation*}
X_a' = X_1 + X_2 + X_3, \qquad X_i = \sum_{e\in \mathcal{C}_i} \frac{\chi(x_{t(e)}y_{h(e)})}{d_{t(e)}d_{h(e)}}.
\end{equation*}
We can do a similar decomposition for $X_b'$.

\medskip\noindent
{\bf 3.}
The following martingale lemma from \cite{10.1145/73007.73063,broderOptimalConstructionEdgeDisjoint1998} is the core estimate that we want to apply to each of the $X_i$'s.

\begin{lemma}
\label{lem:martingale}
Let be $G_1$ and $G_2$ two graphs generated via the configuration model through perfect matchings $P_1$ and $P_2$. Let $\{e^{(1)}_i\}_{i \ge 1}$ be the edges of $G_1$ and $\{e^{(2)}_i\}_{i\ge 1}$ be the edges of $G_2$, ordered as above. Define an equivalence relation on the probability space $\Omega$ by setting $G_1 \equiv_k G_2$ when $\{e^{(1)}_i\}_{i=1}^{k}=\{e^{(2)}_i\}_{i=1}^{k}$, i.e., the first $k$ pairings match. Let $\Omega_k$ be the set of equivalence classes, and $\mathcal{F}_k$ the corresponding $\sigma$-algebra with $\mathcal{F}_0=\mathcal{F}$. Consider a bounded measurable $f\colon\,\mathcal{G}(\vec{d}) \to \R$, and set $Y_k = \Exp[f \mid \mathcal{F}_k]$. Note that 
\begin{itemize}
\item[]
$(Y_k)_{0\le k\le m_1/2}$ is a Doob martingale with $\Exp[Y_k \mid \mathcal{F}_{k-1}] =Y_{k-1}$\\ 
and $Y_0=\Exp[f]$, $Y_{m_1/2}=f$. 
\end{itemize}
Define $Z_k=Y_k-Y_{k-1}$, and suppose that there exist functions $(g_k(\zeta))_{1\le k\le \tfrac12 m_1}$ such that 
\begin{equation*}
\Exp\left[\eee^{\zeta ^2 Z^2_k} \mid \mathcal{F}_{k-1}\right] \leq g_k(\zeta), \qquad 1 \leq k \leq m_1/2.
\end{equation*}
Then, for all $t \geq 0$ and $\zeta>0$,
\begin{equation}
\label{eq:martingaleconc}
\Pro\lep\left|f-\Exp[f]\right|\geq t\rip\leq 2\,\eee^{-\zeta t}\prod_{k=1}^{m_1/2} g_k(\zeta).
\end{equation}
\end{lemma}

\begin{remark}
{\rm The condition on the existence $g_k(\zeta)$ can be rephrased as the existence of a random variable $W_k$ that stochastically dominates $Z_k$ on $(\Omega_{k-1},\mathcal{F}_{k-1})$ (see \cite{10.1145/73007.73063}).}\hfill$\spadesuit$
\end{remark}

\begin{lemma}{\cite[Lemma 15]{broderOptimalConstructionEdgeDisjoint1998}}\label{lem:concentration}
There exist constants $B_\ell>0$, $\ell = 1,2,3$, depending only on the ratio $\frac{m_\infty}{m_0}$, such that
\begin{equation}
\label{eq:conc}
\Pro\lep|X_\ell-\Exp[X_\ell]| \geq \frac{t}{m_\infty^{3/2}}\rip \leq 2\,\eee^{-tn+B_\ell n}, 
\qquad \ell = 1,2,3.
\end{equation}
\end{lemma}

\begin{proof}
While the properties in $\mathcal{C}_3$ allow us to apply standard martingale arguments to capture $X_3$ (see below), the properties in $\mathcal{C}_1$ and $\mathcal{C}_2$ force us to use Lemma ~\ref{lem:martingale} to capture $X_1$ and $X_2$. 

\medskip\noindent
{\bf i.}
From the definition of $\mathcal{C}_1$ and $\mathcal{C}_2$ it follows that a bound on $X_1$ implies by symmetry a bound on $X_2$ (with, possibly, different constants). We will therefore focus on $X_1$, the result for $X_2$ carrying through trivially. Without loss of generality we may reorder the indices in such a way that $|x_i| \geq |x_{i+1}|$ and for each $e_i=\{\alpha_i, \beta_i\}$ use the lexicographic order $\alpha_{i+1} > \alpha_i$ and $\beta_i > \alpha_i$ (i.e., we perform the pairing sequentially from left to right; see \cite[Lemma 7.6]{vanderhofstadRandomGraphsComplex2016}). It follows that $v(\alpha_i) \leq v(\alpha_{i+1})$ and $|x_{v(\alpha_i)}| \geq |x_{v(\alpha_{i+1})}|$. Define 
\begin{equation*}
\hat{\chi}(x,y) = \begin{cases}
xy, \quad &\text{if } |xy|< \frac{\sqrt{m_\infty}}{n}\text{ and } |x|>\frac{1}{\epsilon \sqrt{n}},\\
0, \quad &\text{otherwise}.
\end{cases}
\end{equation*}
For each $e=\{\alpha,\beta\}$ in the configuration $\mathcal{C}$ we have
\begin{equation*}
X_1=\sum_{e\in \mathcal{C}}q(e)
\end{equation*}
with $q(e)=\hat{\chi}(x_{v(\alpha)},y_{v(\beta)})/d_{v(\alpha)}d_{v(\beta)}$.

\medskip\noindent
{\bf ii.}
Next, take $Y_k = \Exp[X_1 \mid \mathcal{F}_k]$, with $Y_0 = \Exp[X_1]$ and $Y_m=X_1$. Then $(Y_k)_{k\in\N_0}$ is a Doob martingale and, by the definition of the configuration model, we can write $Z_k=Y_k-Y_{k-1}$ as 
\begin{equation*}
Z_k(\mathcal{C}) = \frac{2^{\tfrac12 m_1-k}(\tfrac12 m_1-k)!}{(m_1-2k)!}\lep 
\sum_{\mathcal{C}'\equiv_k \mathcal{C}}X_1(\mathcal{C}')
-\frac{1}{m_1 -2k +1}\sum_{\mathcal{C}''\equiv_{k-1}\mathcal{C}} X_1(\mathcal{C}'')\rip.
\end{equation*}
Now that we have an expression for $Z_k$, we can use the method of switching (see, for example, \cite{wormaldModelsRandomRegular1999}). Indeed, given a $\mathcal{C}' \equiv_k \mathcal{C}$, we can define a quantity $\mathcal{C}'_\eta$ as follows. Given the first $k$ pairings of $\mathcal{C}$, let $I$ be the set of points already paired, and let $\{\alpha,\beta\}$ be the $k$-th pair. Put $\eta \notin I-\{\beta\}$ and $\{\eta, \gamma\} \in \mathcal{C}'$. Then $\mathcal{C}'_\eta$ is the pairing obtained from $\mathcal{C}'$ by mapping 
\begin{equation*}
\{\alpha,\beta\},\{\eta, \gamma\} \to \{\alpha,\eta\},\{ \gamma,\beta\}. 
\end{equation*}
Is easy to see that $\mathcal{C}'_\eta\equiv_{k-1} \mathcal{C}$ and that $\{\{\mathcal{C}'_\eta\colon\,\eta \notin I-\{\beta\}\} \mid \mathcal{C}'\equiv_{k} \mathcal{C}\}$ is a partition of $\{\mathcal{C}'' \mid \mathcal{C}''\equiv_{k-1}\mathcal{C}\}$. We can therefore rewrite 
\begin{equation*}
\begin{split}
Z_k(\mathcal{C})
&=\frac{2^{\tfrac12 m_1-k}(\tfrac12 m_1-k)!}{(m_1-2k)!} \sum_{\mathcal{C}'\equiv_k \mathcal{C}}
\sum_{\eta\notin I} (X_1(\mathcal{C}')-X_1(\mathcal{C}'_\eta))\\
&=\sum_{\eta\notin I} \sum_{\gamma\notin I, \gamma\neq\eta} 
\frac{q(\{\alpha,\beta\}) + q(\{\eta, \gamma\}) - q(\{\alpha,\eta\}) - q(\{\gamma,\beta\})}{(2m-2k+1)(2m-2k-1)}.
\end{split}
\end{equation*}
Because $\sum_i x_u^2 \leq 1$, there are at most $\varepsilon^2 n$ indices of $x$ such that $|x_i| > 1/(\varepsilon \sqrt{n})$. By the definition of $X_1$, the lexicographic ordering of $\{\alpha_i,\beta_i\}$ and the ordering of $|x_j|>|x_{j+1}|$, there exists a $\tilde{k}$ such that $Z_{\tilde{k}}=0$. Take $\tilde{k}=\varepsilon^2m_\infty n$. For $k>\tilde{k}$, we have that 
\begin{equation*}
m_1- 2k -1 \geq m_1 - 2\varepsilon^2 m_\infty n - 1 \geq m_0 n,
\end{equation*}
where the free parameter $\varepsilon$ has to be fixed such that the last inequality holds. (Note that, because there exists a constant $\theta$ such that $\frac{m_\infty}{m_0} < \theta$, we can always choose an $\varepsilon$ small enough so that this holds). Hence we can bound
\begin{equation*}
|Z_k(\mathcal{C})|\leq\frac{1}{(m_0n)^2} \sum_{\eta \notin I} \sum_{\gamma\notin I , \gamma \neq \eta} \lep 
|q(\{\alpha,\beta\})| + |q(\{\eta, \gamma\})| + |q( \{\alpha,\eta\})| + |q(\{ \gamma,\beta\})| \rip.
\end{equation*}

\medskip\noindent
{\bf iii.}
Define 
\begin{equation*}
y^\alpha=\frac{1}{|x_{v(\alpha)}|}\text{min}\left\{ |yx_{v(\alpha)}|,\frac{\sqrt{m_\infty}}{n}\right\}.
\end{equation*}
Because $\alpha<\beta$, we can bound
\begin{equation*}
q(\{\alpha,\beta\})=\frac{\hat{\chi}(x_{v(\alpha)},y_{v(\beta)})}{d_{v(\alpha)}d_{v(\beta)}}
\leq \frac{y^\alpha_{v_\beta}|x_{v(\alpha)}|}{\dmin^2}.
\end{equation*}
A similar bound holds for $\{\alpha,\eta\}$ for the same reason. For the other two edges, $\{\eta, \gamma\}$ and $\{ \gamma,\beta\}$, we need to upper bound with a symmetric term, because we do not know whether $\gamma>\eta$ or $\gamma<\eta$. Thus, we have the upper bound
\begin{equation*}
q(\{\eta, \gamma\}) \leq \frac{1}{\dmin^2}\lep y^\gamma_{v_\eta} |x_{v(\gamma)}| + y^\eta_{v_\gamma}| x_{v(\eta)}|\rip
\end{equation*} 
(the same bound holds for $\{\gamma,\beta\}$). Moreover, by the lexicographic order, $x_{v(\alpha)}$ bounds all the other components, and therefore $y^\alpha_{v_\beta}|x_{v(\alpha)} \leq y^\alpha_{v_\beta}|x_{v(\alpha)}$.  Now note that $\sum_i |y_i|\leq \sqrt{n}$ (because $\sum_i y_i^2 \leq 1$), and so
\begin{equation*}
\sum_{\eta\notin I} y^\alpha_{v(\eta)} \leq \sum_{\eta\notin I} |y_{v(\eta)}| \leq m_\infty \sum_i |y_i|\leq m_\infty \sqrt{n}.
\end{equation*}
By the previous considerations, substituting into the expression for $Z(\mathcal{C})$, we have
\begin{equation*}
\begin{split}
|Z_k(\mathcal{C})|
&\leq\frac{1}{m_0^2}\lep y^\alpha_{v_\beta} |x_{v(\alpha)}| + \lep y^\alpha_{v_\eta}|x_{v(\alpha)}| + y^\alpha_{v_\gamma}| 
x_{v(\alpha)}|\rip + y^\alpha_{v_\eta} |x_{v(\alpha)}| + \lep y^\alpha{v_\beta} |x_{v(\alpha)}|
+ y^\alpha_{v_\gamma}| x_{v(\alpha)}|\rip\rip\\
&\leq \frac{4 m_\infty^2}{m_0^4} |x_{v(\alpha)}|\lep y^\alpha_{v(\beta)} + \frac{1}{\sqrt{n}}\rip.
\end{split}
\end{equation*}

\medskip\noindent
{\bf iv.}
From the above bounds we are able to obtain an upper bound for $\Exp[\exp(\zeta^2 Z_k^2) \mid \mathcal{F}_{k-1}]$ and then use lemma \ref{lem:martingale}. Indeed,
\begin{equation*}
\Exp\left[\eee^{\zeta^2 Z_k^2} \mid \mathcal{F}_{k-1}\right]
\leq \frac{1}{m_1 -2k -1} \sum_{\omega \notin I \backslash \beta} 
\exp\left[16 \zeta^2 m_\infty^4 m_0^{-8} (x_{v(\alpha)})^2\right]
\lep y^\alpha_{v(\omega)} + \frac{1}{\sqrt{n}}\rip^2.
\end{equation*}
Looking at \eqref{eq:martingaleconc}, we see that we have to fix $\zeta = m_\infty^{3/2} n$ in order to achieve the required bound. Using that $x_{v(\alpha)} \leq \sqrt{m_\infty}/(\epsilon \sqrt{n})$ (because $x_{v(\alpha)} y^\alpha_{v(\beta)} \leq \sqrt{m_\infty}{n}$ and $y^\alpha_{v(\beta)}\geq \varepsilon/ \sqrt{n}$), the exponent in the previous display is bounded by $64 \theta^8 \varepsilon^{-2}$. Using that $\eee^x \leq 1 + x\eee^x$, $x\geq 0$, and putting $B = 64\theta^8\epsilon^{-2}$, we have
\begin{equation*}
\begin{split}
\Exp\left[e^{\zeta^2 Z_k^2}|\mathcal{F}_{k-1}\right]&\leq 1+\frac{B}{m_1 -2k -1} 
\sum_{\omega \notin I \backslash \beta } \zeta^2 m_\infty ^4  m_0^{-8} (x_{v(\alpha)})^2\lep 
y^\alpha_{v(\omega)} + \frac{1}{\sqrt{n}}\rip^2\\
&\leq 1+B \zeta^2 m_\infty ^4  m_0^{-9} (x_{v(\alpha)})^2 \sum_{\omega }\lep y^2_{v(\omega)}
+2 y_{v(\omega)} \frac{1}{\sqrt{n}}+\frac{1}{n}\rip\\
&\leq \exp\left[4 B \frac{\zeta^2 m_\infty^5} {m_0^9n}(x_{v(\alpha)})^2\right].
\end{split}
\end{equation*}

Next, let us pick an index $i(k)$ such that, for all $\mathcal{C}\in \Omega$,
\begin{equation*}
\Exp\left[e^{\zeta^2 Z_k^2} \mid \mathcal{F}_{k-1}\right] \leq 
\exp\left[4 B \frac{\zeta^2 m_\infty^5}{m_0^9n}(x_{i(k)})^2\right].
\end{equation*}
One possible choice is to take $i(k) = \lceil k/m_\infty \rceil$. We finally get
\begin{equation*}
\Pro\lep|X_1-\Exp[X_1]| \geq \frac{t}{m_\infty^{3/2}}\rip 
\leq 2\,\eee^{-\frac{\zeta t}{m_\infty^{3/2}}} \eee^{\sum_{k=1}^{\tfrac12 m_1} 
4 B \frac{\zeta^2 m_\infty^5} {m_0^9n}(x_{i(k)})^2}
\leq 2\, \eee^{-tn + 4B\frac{m_\infty^9}{m_0^9}n},
\end{equation*}
which proves what was said at the beginning of the proof \eqref{eq:conc} for $\ell=1,2$.

\medskip\noindent
{\bf v.}
Finally, consider $X_3$. In view of the bounds in $\mathcal{C}_3$, this case can be dealt with via classical martingale arguments (see, for example, the McDiarmid inequality and its generalizations in \cite{boucheronConcentrationInequalitiesNonasymptotic2013b}). Considering the variables $Y_k=\Exp[X_3 \mid \mathcal{C}_3]$, we have that $|Y_k-Y_{k-1}|\leq 4/(\varepsilon^2 n m_0^2)$. Thus, given our choice of $\varepsilon=\varepsilon(\theta)$ being constant, we have
\begin{equation*}
\Pro\lep |X_3-\Exp[X_3]|\geq t m_{\infty}^{-3/2}\rip \leq 2\,\eee^{-\frac{t^2 \varepsilon^4n^2m_0^4}{16 m_{\infty}^3 m_1}}
\leq 2\,\eee^{-C(\varepsilon,\theta)t^2n},
\end{equation*}
where $C(\varepsilon,\theta)>0$.
\end{proof}

\medskip\noindent
{\bf 4.}
Combining the results for $\ell=1,2,3$ in the above lemma, we get an exponential bound on $X'_a$ of the type
\begin{equation*}
\Pro\lep|X'_a-\Exp[X'_a]| \geq \frac{t}{m_\infty^{3/2}}\rip \leq 2\eee^{-tn+B_a n},
\end{equation*}
where $B_a$ is a suitable constant, possibly different from any of the $B_\ell$. By symmetry, the same holds for $X'_b$, which proves the concentration around $m_\infty^{-3/2}$ of $X'$. 

\begin{remark} 
{\rm Up to now we have worked with multi-graphs, so we have to pass to simple graphs with a prescribed degree sequence. Is easy to see that, because of \eqref{lem:simpleCM_1}, we can express the final result as saying that, conditional on the event that the graph is simple, there exist two constants $\hat{\xi}<\xi$ in $(0,1)$ and a constant $K(\theta,\xi)>0$ such that
\begin{equation*}
\Pro\lep|X'-\Exp[X']| \geq K{m_\infty^{-3/2}}\rip \leq \eee^{O(\hat{d}^2)} \hat{\xi}^n \leq \xi^n.
\end{equation*}
}\hfill$\spadesuit$
\end{remark}


\subsection{Estimate of second contribution}
\label{sec:second}

It remains to show that the pairs with $x,y\notin \mathcal{B}$ give a bounded contribution of the order $O(m_\infty^{-3/2})$ with sufficiently high probability. For this it suffices to show that a simple random graph $G$ with a prescribed degree sequence $\vec{d}$ cannot have too dense subgraphs (see \cite[Lemma 2.5]{10.1145/73007.73063} and \cite[Lemma 16]{broderOptimalConstructionEdgeDisjoint1998}):

\begin{lemma}
\label{lem:szemeredi}
Let $G$ be a simple random graph of size $n$ drawn uniformly at random with a given degree sequence $\vec{d}$. Let $A,B \subseteq [n]$ be two subsets of the vertex set, and let $e(A,B)$ be the set of edges $e = \{\alpha,\beta\}$ such that either $\alpha \in A, \beta \in B$ or $\alpha \in B, \beta \in A$. Since $\mu(A,B) = \theta |A||B| \frac{m_\infty}{n}$ with $\theta>m_\infty/m_0$ a sufficiently large constant, for any $K>0$ there exist a constant $C=C(\theta,K)$ such that with probability $1-o(n^{-K})$ any pair $A,B$ with $|A|\leq |B|$ satisfies at least one of the following:
\begin{align}
\label{eq:szemeredi_1}
&e(A,B)\leq C\mu(A,B)\\ 
\label{eq:szemeredi_2}
&e(A,B)\log\lep\frac{e(A,B)}{\mu(A,B)}\rip \leq C |B| \log \lep\frac{n}{|B|}\rip.
\end{align}
\end{lemma}  

The above lemma has the following corollary.

\begin{corollary}
For all $x,y\in T$,
\begin{equation*}
X''=O\left(\frac{1}{m_\infty^{3/2}}\right) \quad \textnormal{with probability at least } 1-O(n^{-K}),
\end{equation*}
where $K$ is the constant in Lemma \ref{lem:szemeredi}.
\end{corollary}

\begin{proof}
Fix $x,y \in T$ and define
\begin{equation*}
\begin{split}
S_i(x) &=\left\{\ell\colon\, \frac{\varepsilon^{2-i}}{\sqrt{n}} \leq |x_\ell| <\frac{\varepsilon^{1-i}}{\sqrt{n}} \right\}, \quad i\in I,\\
S_j(y) &=\left\{\ell\colon\, \frac{\varepsilon^{2-i}}{\sqrt{n}} \leq |y_\ell| <\frac{\varepsilon^{1-i}}{\sqrt{n}} \right\}, \quad j\in J,
\end{split}
\end{equation*}
where $I=\{i|S_i(x)\neq \emptyset\}$ ($J$ is defined similarly). Then, for $x\in T$, write
\begin{equation*}
x_u|_S = \begin{cases}
x_u, \quad \textnormal{if } u\in S,\\
0,  \quad \textnormal{otherwise.}
\end{cases}
\end{equation*}
In order to apply Lemma \ref{lem:szemeredi}, define $A_i = S_i(x)$ and $B_j = S_j(y)$, and let $a_i$ and $b_i$ be their cardinality, respectively. Divide the set of indices into two groups:
\begin{equation*}
\begin{aligned}
\mathcal{E} &= \big\{(i,j) \colon\, i,j>0, \varepsilon^{2-i-j}>\sqrt{m_\infty}, a_i\leq b_j\big\},\\
\mathcal{E}' &= \big\{(i,j) \colon\,i,j>0, \varepsilon^{2-i-j}>\sqrt{m_\infty}, a_i> b_j\big\}.
\end{aligned}
\end{equation*}
By the definition of $X''$ and the set $\mathcal{B}$, we have
\begin{equation*}
X'' = \sum_{x_i y_j>\sqrt{m_\infty}/n} x_i A_{ij} y_j
= \sum_{(i,j) \in \mathcal{E}} (x|_{A_i})^tAy|_{B_j} + \sum_{(i,j)\in \mathcal{E}'} (x|_{A_i})^tAy|_{B_j}.
\end{equation*}
It suffices to show that either of the contributions coming from $\mathcal{E}$ or $\mathcal{E}'$ are $O(m_\infty^{-3/2})$ (the other will follow by symmetry). Focus on $\mathcal{E}$. Putting $e_{i,j}=e(A_i,B_j)$ and $\mu_{i,j}=\mu(A_i,B_j)$, we see that the bound can be rewritten as
\begin{equation*}
\frac{1}{n} \sum_{(i,j) \in \mathcal{E}} \frac{e_{i,j}}{\varepsilon^{i+j}} = O(\sqrt{m_\infty}\,).
\end{equation*}
Divide $\mathcal{E}$ into the union of $\mathcal{E}_a$ and $\mathcal{E}_b$, where $\mathcal{E}_a$ satisfies \eqref{eq:szemeredi_1} and $\mathcal{E}_b$ satisfies \eqref{eq:szemeredi_2}. Clearly $\mathcal{E}=\mathcal{E}_a\cup\mathcal{E}_b$ and 
\begin{equation*}
\frac{1}{n} \sum_{(i,j) \in \mathcal{E}} \frac{e_{i,j}}{\varepsilon^{i+j}}\leq \frac{1}{n} \sum_{(i,j) \in \mathcal{E}_a} 
\frac{e_{i,j}}{\varepsilon^{i+j}}+\frac{1}{n} \sum_{(i,j) \in \mathcal{E}_b} \frac{e_{i,j}}{\varepsilon^{i+j}}.
\end{equation*}
If we are able to show that both contributions from $\mathcal{E}_a$ and $\mathcal{E}_b$ are $O(\sqrt{m_\infty}\,)$, then the theorem follows.
It is easy to see that $\mathcal{E}_a$ gives a bounded contribution. Indeed,
\begin{equation*}
\frac{1}{n} \sum_{(i,j) \in \mathcal{E}_a} \frac{e_{i,j}}{\varepsilon^{i+j}} 
\leq \frac{1}{n^2} \sum_{(i,j)\in\mathcal{E}_a} \frac{Ca_i b_j \theta m_\infty}{\varepsilon^{i+j}}
\leq \frac{C' m_\infty}{n^2} \sum_{(i,j)\in\mathcal{E}_a} \frac{a_i b_j}{\varepsilon^{2(i+j)}\sqrt{m_\infty}}
= O(\sqrt{m_\infty}\,),
\end{equation*}
where in the last step we use that, because $\sum_i x_i^2\leq 1$,
\begin{equation*}
\sum_{i\in I} \frac{a_i}{\varepsilon^{2(i-2)}} \leq n, \qquad \sum_{j\in J} \frac{b_j}{\varepsilon^{2(i-2)}} \leq n.
\end{equation*}
It remains to show that
\begin{equation*}
\frac{1}{n} \sum_{(i,j)\in \mathcal{E}_b} \frac{e_{i,j}}{\varepsilon^{i+j}}=O(\sqrt{m_\infty}\,).
\end{equation*}
In order to do so, we divide $\mathcal{E}_b$ into subsets $\mathcal{E}_b^{(\ell)}$, $\ell=1,\dots,5$, having the following properties:
\begin{equation*}
\begin{split}
&(1)\quad \varepsilon^j\geq \varepsilon^i\sqrt{m_\infty}.\\
&(2)\quad e_{i,j}\leq\frac{ \mu_{i,j}}{\varepsilon^{i+j}\sqrt{m_\infty}}.\\
&(3)\quad \log\lep\frac{e_{i,j}}{\mu_{i,j}}\rip\geq \frac{1}{4} \log\lep\frac{n}{b_j}\rip.\\
&(4)\quad \frac{n}{b_j}\leq \eee^{-4j}.\\
&(5)\quad \frac{n}{b_j}> \eee^{-4j}.
\end{split}
\end{equation*}
For $j>i$ we have that $\mathcal{E}_b^{(\ell)}\nsubseteq \mathcal{E}_b^{(\ell)}$ and $\mathcal{E}_b=\cup_\ell \mathcal{E}_b^{(\ell)}$. Thus it suffices to show a bound of $O(\sqrt{n})$ for each of the quantities $S_\ell=1/n \sum_{\mathcal{E}_b^{(\ell)}} e_{i,j}/\varepsilon^{i+j} $.
 
\medskip\noindent 
$\bullet$ For $\ell=1$ we note that, since $e_{i,j}\leq a_i m_\infty$,
\begin{equation*}
S_1\leq \frac{1}{n} \sum_i \sum_{j| \varepsilon^j\geq \varepsilon^i \sqrt{m_\infty}} \frac{a_i m_\infty}{\varepsilon^{i+j}}
= bO\left(\frac{1}{n}\sum_i \frac{a_i\sqrt{m_\infty}}{\varepsilon^{2i}}\right) = O(\sqrt{m_\infty}\,).
\end{equation*}

\medskip\noindent
$\bullet$ For $\ell=2$ we obtain
\begin{equation*}
S_2\leq \frac{1}{n}\sum_{ij}\frac{\mu_{ij}}{\varepsilon^{2(i+j)}\sqrt{m_\infty}}
=O(\frac{\sqrt{m_\infty}}{n^2}\sum_{ij}\frac{a_i b_j}{\varepsilon^{2(i+j)}}) = O(\sqrt{n}).
\end{equation*}

\medskip\noindent
$\bullet$ For $\ell=3$, because the pairs $(i,j)\in \mathcal{E}_b$ have property \eqref{eq:szemeredi_2}, it follows easily that $e_{ij}=O(b_j)$. Furthermore, because $\mathcal{E}^{(3)}_b\nsubseteq \mathcal{E}_b^{(1)}$, we have that $\forall (i,j) \in \mathcal{E}_b^{(3)}$, $e^j<e_i\sqrt{m_\infty}$. It follows that
\begin{equation*}
S_3 = O\left(\frac{1}{n}\sum_j \sum_{i|\varepsilon^i>e^j/\sqrt{m_\infty}} \frac{b_j}{\varepsilon^{i+j}}\right)
= O\left(\frac{1}{n}\sum_j \frac{\sqrt{m_\infty}\,b_j}{\varepsilon^{2j}}\right) = O(\sqrt{m_\infty}\,).
\end{equation*}

\medskip\noindent
$\bullet$ For $\ell=4$ we take advantage of the fact that $(i,j)\in \mathcal{E}_b^{(4)}$ do not belong to $\mathcal{E}_b^{(3)}$ and $\mathcal{E}_b^{(2)}$. This implies that
\begin{equation*}
\frac{e_{ij}}{m_{ij}}\leq\frac{1}{e^j} \qquad \frac{e_{ij}}{\mu_{ij}}\geq \frac{1}{\varepsilon^{i+j}\sqrt{m_\infty}}.
\end{equation*}
Hence we have $\varepsilon^{-i}\leq \sqrt{m_\infty}$ and, by \eqref{eq:szemeredi_2}, also $e_{ij}=O(jb_j)$. We can therefore conclude that
\begin{equation*}
S_4= O\left(\frac{1}{n}\sum_j \sum_{i|\varepsilon^{-i}\leq \sqrt{m_\infty}}\frac{jb_j}{\varepsilon^{i+j}}\right) 
= O\left(\frac{\sqrt{m_\infty}}{n}\sum_j\frac{jb_j}{\varepsilon^j}\right) = O(\sqrt{m_\infty}\,),
\end{equation*}
where in the last equality we use that $\sum_{j\in J}b_j/(n\varepsilon^2)=O(1)$. 

\medskip\noindent
$\bullet$ For $\ell=5$, using the property in (5) and \eqref{eq:szemeredi_2}, we have
\begin{equation*}
e_{ij}\leq Cn\varepsilon^{4j} \log \varepsilon^{-4j} = O(nj\varepsilon^{4j}).
\end{equation*}
Also, using that $ \mathcal{E}_b^{(4)}\nsubseteq \mathcal{E}_b^{(4)}$, we have $\varepsilon^j< \varepsilon ^i\sqrt{m_\infty}$, from which we conclude that
\begin{equation*}
S_5 = O\left(\sum_j\sum_{i|\varepsilon^i>\varepsilon^j/\sqrt{m_\infty}}j \varepsilon^{3j-i}\right)
= O\left(\sqrt{m_\infty}\sum_j j\varepsilon^{2j}\right) = O(\sqrt{m_\infty}\,).
\end{equation*}
This completes the proof.
\end{proof}


\section{Proof of the main theorem}
\label{sec:expansion}


\paragraph{\bf Expansion.}

Throughout this section we abbreviate $A=A_{G_n}$ and condition on the event that $G_n$ is simple. Recall Lemma~\ref{lem:shift}. Compute
\begin{equation*}
\begin{split}
&Av_1=\lambda_1v_1,\\
&\lep H +|\tilde{e}\rangle\langle \tilde{e}| + \lep \Exp[A]-|\tilde{e}\rangle\langle \tilde{e}| \rip\rip v_1=\lambda_1v_1,\\
&\lep H +|\tilde{e}\rangle\langle \tilde{e}| 
-\textnormal{diag}\lep \frac{d_1}{m_1-1},\dots, \frac{d_n}{m_1-1}\rip\rip v_1=\lambda_1v_1.
\end{split}
\end{equation*}
Rewriting the equation we have,
\begin{equation*}
\begin{split}
&\langle \tilde{e}, v_1\rangle \tilde{e} = \lep\lambda_1\mathbb{I}
+ \textnormal{diag}\lep \frac{d_1}{m_1-1},\dots, \frac{d_n}{m_1-1}\rip-H\rip v_1.
\end{split}
\end{equation*}
Therefore componentwise we have the following inequality,
\begin{equation*}
\begin{split}
&\lep\lambda_1+\frac{m_0}{m_1-1}\rip\lep\mathbb{I}-\frac{H}{\lambda_1
+\frac{m_0}{m_1-1}}\rip v_1 \leq \langle \tilde{e}, v_1\rangle \tilde{e} 
\leq \lep\lambda_1+\frac{m_\infty}{m_1-1}\rip\lep\mathbb{I}-\frac{H}{\lambda_1+\frac{m_\infty}{m_1-1}}\rip v_1,
\end{split}
\end{equation*}
If $x$, $y$ and $e$ are non-negative vector (the non-negativity of $v_1$ follows from Perron-Frobenius theory) with $x\geq y$, then $\langle e, x  \rangle\geq \langle e,y \rangle$, we can use Corollary \ref{cor:concentration} and invert the matrix multiplying $v_1$. Indeed, given that $\|H\|=O(\sqrt{m_\infty}\,)$ and $\lambda_1 \sim m_2/m_1$ on the event with probability $1-o(n^{-K})$ of Corollary \ref{cor:concentration} and Theorem \ref{thm:edge}, we can invert and expand
\begin{equation*}
\lep\mathbb{I}-\frac{H}{\lambda_1+\frac{\dmin}{m_1-1}}\rip^{-1}
= \sum_{k\in\N_0} \frac{H^k}{(\lambda_1+\frac{\dmin}{m_1-1})^k},
\end{equation*}
and similarly for $\dmax$. Thus,
\begin{equation*}
\begin{split}
&\lambda_1\leq \frac{m_2}{m_1-1} -\frac{\dmin}{m_1-1}
+ \sum_{k\in\N} \frac{\langle \tilde{e}, H^k \tilde{e}\rangle}{(\lambda_1+\frac{\dmin}{m_1-1})^k}\\
&\lambda_1\geq  \frac{m_2}{m_1-1} -  \frac{\dmax}{m_1-1} 
+ \sum_{k\in\N} \frac{\langle \tilde{e}, H^k \tilde{e}\rangle}{(\lambda_1+\frac{\dmax}{m_1-1})^k}.
\end{split}
\end{equation*}
Our final goal is to determine the expectation $\Exp[\lambda_1]$, which splits as
\begin{equation}
\label{eq:conditioning}
\Exp[\lambda_1] = \Exp[\lambda_1|\mathcal{E}]\,\Pro(\mathcal{E}) + \Exp[\lambda_1|\mathcal{E}^c]\,\Pro(\mathcal{E}^c).
\end{equation}
The event $\mathcal{E}^c$ has probability at most $n^{-K}$, where $K$ is a large arbitrary constant. Thus, given the deterministic bound $\lambda_1 \leq n$, we may focus on $\Exp[\lambda_1|\mathcal{E}]\,\Pro(\mathcal{E})$. In order to do this, we need to be able to handle terms of the type
\begin{equation*}
\begin{aligned}
&\frac{m_2}{m_1-1} + \frac{\Exp[\langle \tilde{e},H \tilde{e} \rangle]}{\frac{m_2}{m_1-1}(1+o(1))}
+\frac{\Exp[\langle \tilde{e}, H^2 \tilde{e} \rangle]}{\lep\frac{m_2}{m_1-1}\rip^2(1+o(1))}
+\sum_{k \in \N \setminus \{1,2\}}\frac{\Exp[\langle \tilde{e},H^k \tilde{e} \rangle]}{\lep\frac{m_2}{m_1-1}\rip^k(1+o(1))}
\end{aligned}
\end{equation*}
Since
\begin{equation*}
\frac{\Exp[\langle \tilde{e}, H^k \tilde{e}\rangle]}{\lep\frac{m_2}{m_1-1}\rip^k(1+o(1))}
\leq\frac{m_\infty^{k/2}}{\lep\frac{m_2}{m_1-1}\rip^{k-1}} = o\left(\frac{1}{\dmin^{k/2-1}}\right),
\end{equation*}
the last sum is $o(1/\sqrt{m_\infty})$, which is an error term. It therefore remains to study $\langle \tilde{e}, H^k \tilde{e} \rangle$, $k=1,2$. The study of these moments for the configuration model is more involved than for random regular graphs. 


\paragraph{\bf Case $k=1$.}

Compute
\begin{equation*}
\begin{split}
\langle \tilde{e}, H \tilde{e}\rangle 
&= \langle \tilde{e}, (A - \Exp[A])\,\tilde{e}\rangle\\
&= \frac{1}{m_1-1} \lep \sum_{ij} d_i d_j a_{ij} - \frac{1}{m_1-1} \sum_{ij} d_i^2 d_j^2\rip
= \frac{\sum_j d_j(\sum_{i\sim j}d_i)}{m_1-1}-\frac{m_2^2}{(m_1-1)^2}.
\end{split}
\end{equation*}
Since
\begin{equation*}
\Exp\left[\sum_j d_j  \sum_{i\sim j} d_i\right]=\sum_{ij} d_j d_i\Exp{a_{ij}}=\frac{m_2^2}{m_1-1}-\frac{m_3}{m_1-1},
\end{equation*}
where the last term comes from the presence of selfloops. It follows that $\Exp[\langle \tilde{e},H^k \tilde{e} \rangle]=O(1/\sqrt{n})$. 


\paragraph{\bf Case $k=2$.}

Compute
\begin{equation*}
\begin{aligned}
\langle \tilde{e}, H^2 \tilde{e}\rangle
&= \langle \tilde{e}, (A - \Exp[A])^2\,\tilde{e}\rangle\\ 
&= \frac{1}{m_1-1}\bigg(\sum_{ijk} d_i d_k a_{ij}a_{jk} - \frac{1}{m_1-1} \sum_{ijk} d_i d^2_k a_{ij}d_j\\
&\qquad -\frac{1}{m_1-1} \sum_{ijk} d^2_i d_jd_k a_{jk} + \frac{1}{(m_1-1)^2} \sum_{ijk} d_i^2 d_j^2 d_k^2\bigg).
\end{aligned}
\end{equation*}
Write
\begin{equation*}
\frac{1}{m_1-1}\sum_{ijk} d_i d^2_k a_{ij}d_j = \frac{m_2}{m_1-1}\sum_{ijk} d_i  a_{ij} d_j 
= \frac{m_2}{m_1-1}\sum_{i} d_i \left(\sum_{i\sim j} d_j\right)
\end{equation*}
(by symmetry the third term is equal) and
\begin{equation*}
\sum_{ijk} d_i d_k a_{ij}a_{jk} = \sum_j \lep \sum_{i\sim j} d_i \rip^2
= \sum_j \sum_{i\sim j} d_i^2 + \sum_k \sum_{\substack{i,j\sim k\\ i\neq j}} d_i d_j
= m_3 +\sum_k \sum_{\substack{i,j\sim k\\ i\neq j}} d_i d_j.
\end{equation*}
Indeed in $\sum_j \sum_{i\sim j} d_i^2$, the summand $d_i^2$ appears exactly $d_i$ times, because the node $i$ has exactly  $d_i$ neighbours, and so $\sum_j \sum_{i\sim j} d_i^2=m_3$. Putting the terms together, we get
\begin{equation*}
\langle \tilde{e}, H^2 \tilde{e} \rangle = \frac{1}{m_1-1}\lep m_3 
+ \sum_k \sum_{\substack{i,j\sim k\\ i\neq j}} d_i d_j 
- 2 \frac{m_2}{m_1-1} \sum_{i} d_i \left(\sum_{i\sim j}  d_j\right) + \frac{m_2^3}{(m_1-1)^2}\rip.
\end{equation*}
Taking expectations, we get
\begin{equation*}
\Exp[\langle \tilde{e},H^2 \tilde{e}\rangle] = \frac{1}{m_1-1}\lep m_3 
+ \Exp\left[\sum_k \sum_{\substack{i,j\sim k\\ i\neq j}} d_i d_j - \frac{m^3_2}{(m_1-1)^2}\rip\right].
\end{equation*}

Note that $\Exp[\sum_k \sum_{i,j\sim k,\,i\neq j} d_i d_j]$ is a sum over the wedges centered at vertex $k$, summed all $k$. We can swap the summation over pairs of vertices, and choose a third neighbour to form a wedge, which gives
\begin{equation*}
\sum_k \sum_{\substack{i,j\sim k\\ i\neq j}} d_i d_j 
= \sum_{\substack{i,j\\i\neq j}} d_id_j \sum_k\mathbbm{1}_{(k\sim i,k\sim j)}.
\end{equation*}
Compute 
\small
\begin{equation*}
\begin{split}
&\Exp\left[\sum_{\substack{i,j\\i\neq j}}d_id_j\sum_k\mathbbm{1}_{(k\sim i,k\sim j)}\right]
=\sum_{\substack{i,j\\i\neq j}}d_id_j\sum_k\Exp\left[\mathbbm{1}_{(k\sim i,k\sim j)}\right]\\
&= \sum_{\substack{i,j\\i\neq j}}d_id_j\sum_k\lep\frac{d_i d_j d_k (d_k-1)}{(m_1-1)(m_1-2)}\mathbbm{1}_{k\neq i,k\neq j}
+\frac{d_i d_j (d_k-1) (d_k-2)}{(m_1-1)(m_1-2)}\mathbbm{1}_{k= i \textnormal{ or }k= j}\rip\\
&= \sum_{\substack{i,j\\i\neq j}}d_id_j\sum_k\lep\frac{d_i d_j d_k (d_k-1)}{(m_1-1)(m_1-2)}
-\frac{2d_id_j(d_k-1)}{(m_1-1)(m_1-2)}\mathbbm{1}_{k= i \textnormal{ or }k= j}\rip\\
&= \frac{1}{(m_1-1)(m_1-2)}\sum_{\substack{i,j\\i\neq j}} d^2_id^2_j\sum_k\lep d_k (d_k-1)
-2(d_k-1)\mathbbm{1}_{k= i \textnormal{ or } k= j}\rip\\
&= \frac{1}{(m_1-1)(m_1-2)}\lep(m_2-m_1)\lep\sum_{i,j}d_i^2d_j^2-\sum_i d_i^4\rip
-2\sum_{\substack{i,j\\i\neq j}}d^2_id^3_j-2\sum_{\substack{i,j\\i\neq j}}d^3_id^2_j
+4\sum_{\substack{i,j\\i\neq j}}d_i^2d_j^2\rip\\
&= \frac{1}{(m_1-1)(m_1-2)}\lep(m_2-m_1)\lep m_2^2-m_4\rip-4\lep\sum_{i,j}d^2_id^3_j
-\sum_i d_i^5\rip+4\lep\sum_{i,j}d_i^2d_j^2-\sum_id_i^4\rip\rip\\
&= \frac{1}{(m_1-1)(m_1-2)}\lep(m_2-m_1)\lep m_2^2-m_4\rip-4\lep m_3m_2 -m_5\rip+4\lep m_2^2-m_4\rip\rip\\
&=\frac{m_2^3-m_2 m_4-m_1 m_2^2 +m_1 m_4-4m_3m_2 +4m_5+4m_2^2-4m_4}{(m_1-1)(m_1-2)}.
\end{split}
\end{equation*}
\normalsize
Hence
\small
\begin{equation*}
\begin{aligned}
&\Exp[\langle \tilde{e}, H^2 \tilde{e}\rangle]\\
&=\frac{1}{m_1-1}\bigg( m_3 +\frac{m_2^3-m_2 m_4-m_1 m_2^2 +m_1 m_4-4m_3m_2 
+4m_5+4m_2^2-4m_4}{(m_1-1)(m_1-2)} - \frac{m^3_2}{(m_1-1)^2}\bigg).
\end{aligned}
\end{equation*}
\normalsize
Given the event $\mathcal{E}$, by \eqref{eq:conditioning} and Corollary \ref{cor:lambdaconc}, we have that $\Exp[\lambda_1]$ concentrates around $m_2/m_1$ with an $O(\sqrt{m_\infty})$ error, and so we can write 
\begin{equation*}
\Exp\left[\frac{\langle \tilde{e}, H^2 \tilde{e}\rangle}{\lambda_1^2} \right]
= \frac{\Exp\left[\langle \tilde{e}, H^2 \tilde{e}\rangle \right]}{(\frac{m_2}{m_1-1})^2}(1+o(1)) + o(1).
\end{equation*}

We run through the various contributions separately (using Assumption \ref{ass:deg}). Noting that $nm_0^k \leq m_k \leq nm_\infty^k$ and that. there are positive constants $c,C$ such that $c \leq \frac{m_\infty}{m_0} \leq C$, we have
\begin{equation*}
\begin{aligned}
&\frac{m_2^3}{m_2^2(m_1-1)(m_1-2)} =\Theta\lep\frac{1}{n}\rip,\\
&\frac{m_2m_4}{m_2^2(m_1-2)} = o\left(\frac{1}{\sqrt{n}}\right), \quad
\frac{m_1m_4}{m_2^2(m_1-2)} = \Theta\lep\frac{1}{n}\rip, \quad
\frac{4m_3m_2}{m_2^2(m_1-2)} = o\left(\frac{1}{n}\right),\\
&\frac{4m_5}{m_2^2(m_1-2)} = o\left(\frac{1}{n^{3/2}}\right), \quad 
\frac{4m_2^2}{m_2^2(m_1-2)} = o\left(\frac{1}{n}\right), \quad
\frac{4m_4}{m_2^2(m_1-2)} = o\left(\frac{1}{n^2}\right).
\end{aligned}
\end{equation*}
Therefore
\begin{equation*}
\frac{\Exp\left[\langle \tilde{e}, H^2 \tilde{e}\rangle \right]}{(m_2/(m_1-1))^2}
=\frac{m_3m_1}{m^2_2}-1+o(\frac{1}{\sqrt{n}}),
\end{equation*}
which settles \eqref{eq:mic}.


\paragraph{\bf Weak law of large numbers.}

We want to show that 
\begin{equation*}
\frac{\lambda_1}{\Exp{\lambda_1}}\to 1
\end{equation*}
in $\Pro$-probability. Using Corollary \ref{cor:concentration} and the Weyl interlacing inequality, we have that, with probability $1-n^{-K}$ for $K>0$,
\begin{equation*}
\frac{m_2}{m_1}-O\left(\sqrt{m_\infty}\,\right) \leq \lambda_1\leq \frac{m_2}{m_1} + O(\sqrt{m_\infty}\,).
\end{equation*}
By \eqref{eq:mic},
\begin{equation*}
\frac{\frac{m_2}{m_1}-O(\sqrt{m_\infty}\,)}{\frac{m_2}{m_1}(1+o(1))}\leq\frac{\lambda_1}{\Exp[\lambda_1]}
\leq \frac{\frac{m_2}{m_1}+O(\sqrt{m_\infty}\,)}{\frac{m_2}{m_1}(1+o(1))}.
\end{equation*}
It follows that 
\begin{equation*}
1-O\left(\frac{1}{\sqrt{m_\infty}}\right))\leq \frac{\lambda_1}{\Exp[\lambda_1]} \leq 1+O\left(\frac{1}{\sqrt{m_\infty}}\right)
\end{equation*}
with probability $1-n^{-K}$, and hence the claim follows.


\bibliography{bib}{}

\begin{thebibliography}{10}

\bibitem{alonEigenvaluesExpanders1986}
N.~Alon.
\newblock Eigenvalues and expanders.
\newblock {\em Combinatorica}, 6(2):83--96, June 1986.

\bibitem{arizmendie.TransformSymmetricProbability2009}
O.~Arizmendi~E. and V.~Pérez-Abreu.
\newblock The \${S}\$-transform of symmetric probability measures with
  unbounded supports.
\newblock {\em Proceedings of the American Mathematical Society},
  137(09):3057--3057, Sept. 2009.

\bibitem{BHKY}
R.~Bauerschmidt, J.~Huang, A.~Knowles, and H.-T. Yau.
\newblock Edge rigidity and universality of random regular graphs of
  intermediate degree.
\newblock {\em Geom. Funct. Anal.}, 30(3):693--769, 2020.

\bibitem{bercoviciFreeConvolutionMeasures1993}
H.~Bercovici and D.~Voiculescu.
\newblock Free {Convolution} of {Measures} with {Unbounded} {Support}.
\newblock {\em Indiana University Mathematics Journal}, 42(3), 1993.

\bibitem{bordenaveNewProofFriedman2019}
C.~Bordenave.
\newblock A new proof of {Friedman}'s second eigenvalue {Theorem} and its
  extension to random lifts, Mar. 2019.
\newblock arXiv:1502.04482 [math].

\bibitem{Bordenave:newproof}
C.~Bordenave.
\newblock A new proof of {F}riedman's second eigenvalue theorem and its
  extension to random lifts.
\newblock {\em Ann. Sci. \'{E}c. Norm. Sup\'{e}r. (4)}, 53(6):1393--1439, 2020.

\bibitem{Bordenave:LeLarge}
C.~Bordenave and M.~Lelarge.
\newblock Resolvent of large random graphs.
\newblock {\em Random Structures Algorithms}, 37(3):332--352, 2010.

\bibitem{boucheronConcentrationInequalitiesNonasymptotic2013b}
S.~Boucheron, G.~Lugosi, and P.~Massart.
\newblock {\em Concentration inequalities: a nonasymptotic theory of
  independence}.
\newblock Oxford University Press, Oxford, 1st ed edition, 2013.
\newblock OCLC: ocn818449985.

\bibitem{broderSecondEigenvalueRandom1987}
A.~Broder and E.~Shamir.
\newblock On the second eigenvalue of random regular graphs.
\newblock In {\em 28th {Annual} {Symposium} on {Foundations} of {Computer}
  {Science} (sfcs 1987)}, pages 286--294, Los Angeles, CA, USA, Oct. 1987.
  IEEE.

\bibitem{broderOptimalConstructionEdgeDisjoint1998}
A.~Z. Broder, A.~M. Frieze, S.~Suen, and E.~Upfal.
\newblock Optimal {Construction} of {Edge}-{Disjoint} {Paths} in {Random}
  {Graphs}.
\newblock {\em SIAM Journal on Computing}, 28(2):541--573, Jan. 1998.

\bibitem{Chakrabarty:Hazra:denHollander:Sfragara}
A.~Chakrabarty, R.~S. Hazra, F.~den Hollander, and M.~Sfragara.
\newblock Spectra of adjacency and {L}aplacian matrices of inhomogeneous
  {E}rdos-{R}enyi random graphs.
\newblock {\em Random Matrices Theory Appl.}, 10(1):Paper No. 2150009, 34,
  2021.

\bibitem{demboEmpiricalSpectralDistributions2021}
A.~Dembo, E.~Lubetzky, and Y.~Zhang.
\newblock Empirical {Spectral} {Distributions} of {Sparse} {Random} {Graphs}.
\newblock In M.~E. Vares, R.~Fernández, L.~R. Fontes, and C.~M. Newman,
  editors, {\em In and {Out} of {Equilibrium} 3: {Celebrating} {Vladas}
  {Sidoravicius}}, volume~77, pages 319--345. Springer International
  Publishing, Cham, 2021.
\newblock Series Title: Progress in Probability.

\bibitem{dionigiSpectralSignatureBreaking2021}
P.~Dionigi, D.~Garlaschelli, F.~d. Hollander, and M.~Mandjes.
\newblock A spectral signature of breaking of ensemble equivalence for
  constrained random graphs.
\newblock {\em Electronic Communications in Probability}, 26(none), Jan. 2021.
\newblock arXiv:2009.05155 [math-ph].

\bibitem{dionigiCentralLimitTheorem2023}
P.~Dionigi, D.~Garlaschelli, R.~Subhra~Hazra, F.~den Hollander, and M.~Mandjes.
\newblock Central limit theorem for the principal eigenvalue and eigenvector of
  {Chung}–{Lu} random graphs.
\newblock {\em Journal of Physics: Complexity}, 4(1):015008, Mar. 2023.

\bibitem{erdosSpectralStatisticsErdos2013b}
L.~Erdős, A.~Knowles, H.-T. Yau, and J.~Yin.
\newblock Spectral statistics of {Erdős}–{Rényi} graphs {I}: {Local}
  semicircle law.
\newblock {\em The Annals of Probability}, 41(3B), May 2013.

\bibitem{friedmanProofAlonSecond2004}
J.~Friedman.
\newblock A proof of {Alon}'s second eigenvalue conjecture and related
  problems, May 2004.
\newblock arXiv:cs/0405020.

\bibitem{Friedman}
J.~Friedman.
\newblock A proof of {A}lon's second eigenvalue conjecture and related
  problems.
\newblock {\em Mem. Amer. Math. Soc.}, 195(910):viii+100, 2008.

\bibitem{10.1145/73007.73063}
J.~Friedman, J.~Kahn, and E.~Szemerédi.
\newblock On the second eigenvalue of random regular graphs.
\newblock In {\em Proceedings of the twenty-first annual {ACM} symposium on
  theory of computing}, {STOC} '89, pages 587--598, New York, NY, USA, 1989.
  Association for Computing Machinery.
\newblock Number of pages: 12 Place: Seattle, Washington, USA.

\bibitem{furedi1981eigenvalues}
Z.~F{\"u}redi and J.~Koml{\'o}s.
\newblock The eigenvalues of random symmetric matrices.
\newblock {\em Combinatorica}, 1:233--241, 1981.

\bibitem{garlaschelliCovarianceStructureBreaking2018}
D.~Garlaschelli, F.~den Hollander, and A.~Roccaverde.
\newblock Covariance {Structure} {Behind} {Breaking} of {Ensemble}
  {Equivalence} in {Random} {Graphs}.
\newblock {\em Journal of Statistical Physics}, 173(3-4):644--662, Nov. 2018.

\bibitem{juhasz1981spectrum}
F.~Juh{\'a}sz.
\newblock On the spectrum of a random graph.
\newblock {\em Algebraic methods in graph theory}, 1, 1981.

\bibitem{mckay1981expected}
B.~D. McKay.
\newblock The expected eigenvalue distribution of a large regular graph.
\newblock {\em Linear Algebra and its applications}, 40:203--216, 1981.

\bibitem{shamirSharpConcentrationChromatic1987}
E.~Shamir and J.~Spencer.
\newblock Sharp concentration of the chromatic number on random graphs {G}(n,
  p).
\newblock {\em Combinatorica}, 7(1):121--129, Mar. 1987.

\bibitem{squartini2015breaking}
T.~Squartini, J.~de~Mol, F.~den Hollander, and D.~Garlaschelli.
\newblock Breaking of ensemble equivalence in networks.
\newblock {\em Physical review letters}, 115(26):268701, 2015.

\bibitem{touchette2015equivalence}
H.~Touchette.
\newblock Equivalence and nonequivalence of ensembles: Thermodynamic,
  macrostate, and measure levels.
\newblock {\em Journal of Statistical Physics}, 159(5):987--1016, 2015.

\bibitem{tran:vu:wang}
L.~V. Tran, V.~H. Vu, and K.~Wang.
\newblock Sparse random graphs: eigenvalues and eigenvectors.
\newblock {\em Random Structures Algorithms}, 42(1):110--134, 2013.

\bibitem{vanderhofstadRandomGraphsComplex2016}
R.~v.~d. van~der Hofstad.
\newblock {\em Random {Graphs} and {Complex} {Networks}}.
\newblock Cambridge University Press, 1 edition, Nov. 2016.

\bibitem{vanderhofstadRandomGraphsComplex2023}
R.~v.~d. van~der Hofstad.
\newblock {\em Random {Graphs} and {Complex} {Networks} vol 2}, volume~2.
\newblock 2023.

\bibitem{voiculescuMultiplicationCertainNonCommuting1987}
D.~Voiculescu.
\newblock Multiplication of {Certain} {Non}-{Commuting} {Random} {Variables}.
\newblock {\em Journal of Operator Theory}, 18(2):223--235, 1987.

\bibitem{Vu:survey}
V.~H. Vu.
\newblock Recent progress in combinatorial random matrix theory.
\newblock {\em Probab. Surv.}, 18:179--200, 2021.

\bibitem{wormaldModelsRandomRegular1999}
N.~C. Wormald.
\newblock Models of {Random} {Regular} {Graphs}.
\newblock In J.~D. Lamb and D.~A. Preece, editors, {\em Surveys in
  {Combinatorics}, 1999}, pages 239--298. Cambridge University Press, 1
  edition, July 1999.

\end{thebibliography}
\bibliographystyle{abbrv}


\end{document}